\documentclass[amstex,12pt, amssymb]{article}

\usepackage{mathtext}
\usepackage[cp1251]{inputenc}
\usepackage[T2A]{fontenc}
\usepackage[dvips]{graphicx}
\usepackage{amsmath}
\usepackage{amssymb}
\usepackage{amsxtra}
\usepackage{latexsym}
\usepackage{ifthen}

\newcounter{lemma}[section]

\newcounter{corollary}[section]

\newcounter{remark}[section]

\newcounter{theorem}[section]

\newcounter{proposition}[section]

\newcounter{example}

\numberwithin{equation}{section}

\begin{document}

\markboth{\centerline{E.~SEVOST'YANOV,
S.~SKVORTSOV}}{\centerline{LOGARITHMIC H\"{O}LDER CONTINUOUS
MAPPINGS...}}

\def\cc{\setcounter{equation}{0}
\setcounter{figure}{0}\setcounter{table}{0}}

\overfullrule=0pt


\author{EVGENY SEVOST'YANOV, SERGEI SKVORTSOV}

\title{
{\bf LOGARITHMIC H\"{O}LDER CONTINUOUS MAPPINGS AND BELTRAMI
EQUATION}}

\date{\today}
\maketitle

\begin{abstract}
The paper is devoted to the study of mappings satisfying the inverse
Poletsky inequality.  We study the local behavior of these mappings,
moreover, we are most interested in the case when the corresponding
majorant is integrable on some set of spheres of positive linear
measure. The most important result is a logarithmic H\"{o}lder
continuity of such mappings at inner points. As a corollary, we
obtained the existence of a continuous $ACL$-solution of the
Beltrami equation, which is logarithmic H\"{o}lder continuous.
\end{abstract}

\bigskip
{\bf 2010 Mathematics Subject Classification: Primary 30C65;
Secondary 31A15, 31B25}

\section{Introduction}
As is known, mappings with bounded distortion satisfy the inequality
\begin{equation}\label{eq2}
M(\Gamma)\leqslant N(f, A)\cdot K\cdot M(f(\Gamma))\,,
\end{equation}
where $M$ is the modulus of families of paths $\Gamma$ in the domain
$D,$ $$N(y, f, A)\,=\,{\rm card}\,\left\{x\in A: f(x)=y\right\}\,,
\qquad N(f, A)\,=\,\sup\limits_{y\in{\Bbb R}^n}\,N(y, f, A)\,,$$
$A$ is a Borel set in $D,$ and $K\geqslant 1$ is some constant
defined as
$$K={\rm ess \sup}\, K_O(x, f)\,.$$
Here we also used the common notation $K_O(x, f)=\Vert
f^{\,\prime}(x)\Vert^n/J(x, f)$ for $J(x, f)\ne 0;$ $K_O(x, f)=1$
for $f^{\,\prime}(x)=0,$ and $K_O(x, f)=\infty$ for
$f^{\,\prime}(x)\ne 0$ and $J(x, f)=0.$ In addition, regarding
relation~(\ref{eq2}), we could point to the following papers, see
e.g., \cite[Theorem~3.2]{MRV$_1$} or \cite[Theorem~6.7.II]{Ri}. Some
similar condition applies to mappings whose have unbounded outher
dilatations. In particular, the relation
\begin{equation}\label{eq1C}
M(\Gamma)\le \int\limits_{f(E)}K_I\left(y, f^{\,-1},
E\right)\cdot\rho_*^n(y)\,dm(y)\,,
\end{equation}
holds for the so-called mappings with finite length distortion,
where $E$ is an arbitrary measurable subset of $D,$ $\Gamma$ is a
family of paths in $E$ and $\rho_*(y)\in {\rm adm \,}f(\Gamma)$
(see, e.g., \cite[Theorem~8.5]{MRSY$_2$}). In this manuscript, the
main object of the study are mappings which satisfy some more
general inequality than~(\ref{eq1C}). More precisely, in
relation~(\ref{eq1C}), we will assume that there is a more abstract
Lebesgue measurable function $Q,$ and all the problems related to
H\"{o}lder continuity, which are studied in this article, we will
try to connect with the behavior of this function. For estimates
of~H\"{o}lder type for many well-known classes of mappings, such as
mappings with bounded distortion and quasiconformal mappings, see,
e.g.,~\cite[Theorem~3.2.II]{LV}, \cite[Theorem~1.1.2]{Re},
\cite[Theorem~18.2, Remark~18.4]{Va}, \cite[Theorem~3.2]{MRV$_1$}
and \cite[Theorem~1.8]{GG}. Regarding the more general H\"{o}lder
logarithmic continuity, we point out articles and
monographs~\cite[Theorems~1.1.V and 2.1.V]{Suv},
\cite[Theorem~7.4]{MRSY$_2$}, \cite[Theorem~3.1]{MRSY$_3$} and
\cite[Theorem~5.11]{RS}, cf. \cite[Theorems 4 and 5]{Cr}. We
emphasize that, under some very special restrictions, such mappings
are even Lipschitz, but this fact is very rare for mappings of such
a general nature (see, e.g., \cite[Theorem~1.1, Lemma~4.2,
Theorem~4.1]{RSS$_2$}).

\medskip
In what follows, $M$ denotes the $n$-modulus of a family of paths,
and the element $dm(x)$ corresponds to a Lebesgue measure in ${\Bbb
R}^n,$ $n\geqslant 2,$ see~\cite{Va}. For the sets $A, B\subset{\Bbb
R}^n$ we set, as usual,
$${\rm diam}\,A=\sup\limits_{x, y\in A}|x-y|\,,\quad {\rm dist}\,(A, B)=\inf\limits_{x\in A,
y\in B}|x-y|\,.$$
Sometimes, instead of ${\rm dist}\,(A, B),$ we also write $d(A, B),$
if a misunderstanding is impossible. Given sets $E$ and $F$ and a
given domain $D$ in $\overline{{\Bbb R}^n}={\Bbb R}^n\cup
\{\infty\},$ we denote by $\Gamma(E, F, D)$ the family of all paths
$\gamma:[0, 1]\rightarrow \overline{{\Bbb R}^n}$ joining $E$ and $F$
in $D,$ that is, $\gamma(0)\in E,$ $\gamma(1)\in F$ and
$\gamma(t)\in D$ for all $t\in (0, 1).$ Everywhere below, unless
otherwise stated, the boundary and the closure of a set are
understood in the sense of an extended Euclidean space
$\overline{{\Bbb R}^n}.$ Let $x_0\in\overline{D},$ $x_0\ne\infty,$
$$B(x_0, r)=\{x\in {\Bbb R}^n: |x-x_0|<r\}\,, \quad {\Bbb B}^n=B(0, 1)\,,$$
$$S(x_0,r) = \{
x\,\in\,{\Bbb R}^n : |x-x_0|=r\}\,, S_i=S(x_0, r_i)\,,\quad
i=1,2\,,$$
\begin{equation*}\label{eq1**}
A=A(x_0, r_1, r_2)=\{ x\,\in\,{\Bbb R}^n : r_1<|x-x_0|<r_2\}\,.
\end{equation*}
Let $f:D\rightarrow{\Bbb R}^n,$ $n\geqslant 2,$ and let $Q:{\Bbb
R}^n\rightarrow [0, \infty]$ be a Lebesgue measurable function such
that $Q(x)\equiv 0$ for $x\in{\Bbb R}^n\setminus f(D).$ Let
$A=A(y_0, r_1, r_2).$ Let $\Gamma_f(y_0, r_1, r_2)$ denotes the
family of all paths $\gamma:[a, b]\rightarrow D$ such that
$f(\gamma)\in \Gamma(S(y_0, r_1), S(y_0, r_2), A(y_0, r_1, r_2)),$
i.e., $f(\gamma(a))\in S(y_0, r_1),$ $f(\gamma(b))\in S(y_0, r_2),$
and $\gamma(t)\in A(y_0, r_1, r_2)$ for any $a<t<b.$ We say that
{\it $f$ satisfies the inverse Poletsky inequality} at $y_0\in f
(D)$ if the relation
\begin{equation}\label{eq2*A}
M(\Gamma_f(y_0, r_1, r_2))\leqslant \int\limits_A Q(y)\cdot \eta^n
(|y-y_0|)\, dm(y)
\end{equation}
holds for any Lebesgue measurable function $\eta:(r_1,
r_2)\rightarrow [0, \infty]$ such that
\begin{equation}\label{eqA2}
\int\limits_{r_1}^{r_2}\eta(r)\,dr\geqslant 1 \,.
\end{equation}
Note that the first author established the openness and discreteness
of mappings in~(\ref{eq2*A}) under certain conditions on the
function $Q,$ see, e.g., \cite{Sev}. In a more general case, the
performance of these properties is not guaranteed. Note also that
the equicontinuity of homeomorphisms with condition~(\ref{eq2*A})
with somewhat less general constraints on the domain and the
corresponding mappings is studied in detail in~\cite{SevSkv}. In
this manuscript, we focus on mappings with branching.

\medskip
We now formulate the main results of this article. To this end, we
recall a few more definitions. A mapping $f:D\rightarrow {\Bbb R}^n$
is called a {\it discrete} if the preimage $\{f^{\,-1}(y)\}$ of each
point $y\,\in\,{\Bbb R}^n$ consist of isolated points, and {\it
open} if the image of any open set $U\subset D $ is an open set
in ${\Bbb R}^n.$ 
In the extended Euclidean space $\overline{{{\Bbb R}}^n}={{\Bbb
R}}^n\cup\{\infty\},$ we use the so-called chordal metric $h$
defined by the equalities
\begin{equation}\label{eq3C}
h(x,y)=\frac{|x-y|}{\sqrt{1+{|x|}^2} \sqrt{1+{|y|}^2}}\,,\quad x\ne
\infty\ne y\,, \quad\,h(x,\infty)=\frac{1}{\sqrt{1+{|x|}^2}}\,.
\end{equation}
For a given set $E\subset\overline{{\Bbb R}^n},$ we set
\begin{equation}\label{eq9C}
h(E):=\sup\limits_{x,y\in E}h(x, y)\,.
\end{equation}
The quantity $h(E)$ in~(\ref{eq9C}) is called the {\it chordal
diameter} of the set $E.$

For given sets $A, B\subset \overline{{\Bbb R}^n},$ we put
$$h(A, B)=\inf\limits_{x\in A, y\in B}h(x, y)\,,$$
where $h$ is a chordal metric defined in~(\ref{eq3C}).

\medskip
For domains $D, D^{\,\prime}\subset {\Bbb R}^n,$ $n\geqslant 2,$ and
a Lebesgue measurable function $Q:{\Bbb R}^n\rightarrow [0, \infty]$
equal to zero outside the domain $D^{\,\prime},$ we define by
$\frak{R}_Q(D, D^{\,\prime})$ the family of all open discrete
mappings $f:D\rightarrow D^{\,\prime}$ such that
relation~(\ref{eq2*A}) holds for each point $y_0\in D^{\,\prime}.$

\medskip
\begin{theorem}\label{th1}
{\sl Let $D$ and $D^{\,\prime}$ be domains in ${\Bbb R}^n,$
$n\geqslant 2,$ and let $D^{\,\prime}$ be a bounded domain. Suppose
that, for each point $y_0\in D^{\,\prime}$ and for every
$0<r_1<r_2<r_0:=\sup\limits_{y\in D^{\,\prime}}|y-y_0|$ there is a
set $E\subset[r_1, r_2]$ of a positive linear Lebesgue measure such
that the function $Q$ is integrable over the spheres $S(y_0, r)$ for
every $r\in E.$ Then the family of mappings $\frak{R}_Q(D,
D^{\,\prime})$ is equicontinuous at each point $x_0\in D.$ }
\end{theorem}

\medskip
Theorem~\ref{th1} generalizes the result
of~\cite[Theorem~1.5]{SevSkv} to the case when mappings may have
branch points, and the corresponding function $Q$ may turn out to be
non-integrable in the considered domain $D.$ In particular,
Theorem~\ref{th1} implies the following obvious
\medskip
\begin{corollary}\label{cor2}
{\sl Assume that under the conditions of Theorem~\ref{th1}
$$\int\limits_{S(y_0, r)}Q(y)\,d\mathcal{H}^{n-1}(y)<~\infty$$ for
almost all $r\in (0, r_0)$ and any $y_0\in D^{\,\prime}.$ Then the
family of mappings $\frak{R}_Q(D, D^{\,\prime})$ is equicontinuous
at each point $x_0\in D.$ }
\end{corollary}

\medskip
For a domain $D\subset {\Bbb R}^n,$ $n\geqslant 2,$ and a Lebesgue
measurable function $Q:{\Bbb R}^n\rightarrow [0, \infty],$ we define
by $\frak{F}_Q(D)$ the family of all open discrete mappings
$f:D\rightarrow {\Bbb R}^n$ such that relation~(\ref{eq2*A}) holds
for each point $y_0\in f(D).$ In the case where the function $Q$
behaves somewhat better, and the domain $D$ is the unit ball, the
following result holds. Note that the original (local) version of
this result was published by us in~\cite[Theorem~1.1]{SSD}.

\medskip
\begin{theorem}\label{th5}
{\sl Let $n\geqslant 2,$ and let $Q\in L^1({\Bbb R}^n).$ Suppose
that $K$ is a compact set in ${\Bbb B}^n.$ Now, the inequality
\begin{equation}\label{eq2C}
|f(x)-f(y)|\leqslant\frac{C_n\cdot (\Vert
Q\Vert_1)^{1/n}}{\log^{1/n}\left(1+\frac{r_0}{2|x-y|}\right)}
\end{equation}
holds for any $x, y\in K$ and $f\in \frak{F}_Q({\Bbb B}^n),$ where
$\Vert Q\Vert_1$ denotes $L^1$-norm of $Q$ in ${\Bbb R}^n,$ $C_n>0$
is some constant depending only on $n,$ and $r_0=d(K,
\partial {\Bbb B}^n).$ }
\end{theorem}

\medskip
\begin{remark}
If the domain $D$ does not have finite boundary points, then we will
consider the number $r_0$ to be any positive. The said agreement
will apply to the entire subsequent text of the manuscript.
\end{remark}

\medskip
A separate case of Theorems~\ref{th1} and~\ref{th5}, and also the
Corollary~\ref{cor2}, is a situation when $f$ is a homeomorphism in
$D.$ In this case, let us denote $g:=f^{\,- 1}$ and remark that
\begin{equation}\label{eq2D}
g(\Gamma(S(y_0, r_1), S(y_0, r_2), f(D)))=\Gamma_f(y_0, r_1, r_2)\,.
\end{equation}
In fact, if $\gamma\in g(\Gamma(S(y_0, r_1), S(y_0, r_2), f(D))),$
then $\gamma:[a, b]\rightarrow {\Bbb R}^n,$ and
$\gamma=g\circ\alpha,$ $\alpha:[a, b]\rightarrow {\Bbb R}^n$ and
$\alpha(a)\in S(y_0, r_1),$ $\alpha(b)\in S(y_0, r_2),$
$\alpha(t)\in f(D)$ for $a\leqslant t\leqslant b.$ Now,
$\gamma(t)\in D$ for $a\leqslant t\leqslant b$ and
$f(\gamma)=\alpha\in \Gamma(S(y_0, r_1), S(y_0, r_2), f(D)),$ i.e.,
$\gamma\in \Gamma_f(y_0, r_1, r_2).$ Thus, $g(\Gamma(S(y_0, r_1),
S(y_0, r_2), f(D)))\subset\Gamma_f(y_0, r_1, r_2).$ The inverse
inclusion is proved similarly.

\medskip For domains $D, D^{\,\prime}\subset {\Bbb R}^n,$ $n\geqslant
2,$ and a Lebesgue measurable function $Q:{\Bbb R}^n\rightarrow [0,
\infty]$ equal to zero outside the domain $D,$ we define by
$\frak{R}^*_Q(D, D^{\,\prime})$ the family of all homeomorphisms $g$
of $D^{\,\prime}$ onto $D$ such that relation
\begin{equation}\label{eq8B}
M(f(\Gamma(S(x_0, r_1), S(x_0, r_2), D)))\leqslant \int\limits_D
Q(x)\cdot \eta^n (|x-x_0|)\, dm(x)
\end{equation}
holds for $f=g^{\,-1},$ each $x_0\in D,$ any
$0<r_1<r_2<d_0=\sup\limits_{x\in D}|x-x_0|$ and any Lebesgue
measurable function $\eta: (r_1,r_2)\rightarrow [0,\infty ]$ with
the condition~(\ref{eqA2}). Taking into account Theorem~\ref{th1},
Theorem~\ref{th5} and Corollary~\ref{cor2} and
relation~(\ref{eq2D}), we obtain the following statement.

\medskip
\begin{corollary}\label{cor1}
{\sl Let $D$ and $D^{\,\prime}$ be domains in ${\Bbb R}^n,$
$n\geqslant 2,$ and let $D$ be a bounded domain. Suppose that
$Q:{\Bbb R}^n\rightarrow [0, \infty]$ is a Lebesgue measurable
function and, besides that, for each point $x_0\in\overline{D}$ and
for every $0<r_1<r_2<r_0:=\sup\limits_{x\in D}|x-x_0|$ there is a
set $E\subset[r_1, r_2]$ of positive linear Lebesgue measure such
that the function $Q$ is integrable over the spheres $S(x_0, r)$ for
every $r\in E.$ Then the family of mappings $\frak{R}^*_Q(D,
D^{\,\prime})$ is equicontinuous at each point $y_0$ of
$D^{\,\prime}.$}
\end{corollary}

\medskip
Given domains $D, D^{\,\prime}\subset {\Bbb R}^n,$ $n\geqslant 2,$
and a Lebesgue measurable function $Q:{\Bbb R}^n\rightarrow [0,
\infty],$ we define by $\frak{F}^*_Q(D, D^{\,\prime})$ the family of
all homeomorphisms $g$ of $D^{\,\prime}$ onto $D$ such that
$f=g^{\,-1}$ satisfies relation~(\ref{eq8B}) for each $x_0\in D,$
any $0<r_1<r_2<d_0=\sup\limits_{x\in D}|x-x_0|$ and any Lebesgue
measurable function $\eta: (r_1,r_2)\rightarrow [0,\infty ]$ with
the condition~(\ref{eqA2}). The following statement holds.

\medskip
\begin{corollary}\label{cor4}
{\sl Let $n\geqslant 2,$ and let $Q\in L^1({\Bbb R}^n).$ Suppose
that $K$ is a compact set in ${\Bbb B}^n.$ Now, the inequality
\begin{equation*}\label{eq2CA}
|g(x)-g(y)|\leqslant\frac{C_n\cdot (\Vert
Q\Vert_1)^{1/n}}{\log^{1/n}\left(1+\frac{r_0}{2|x-y|}\right)}
\end{equation*}
holds for every $x, y\in K$ and every $g\in \frak{F}^*_Q(D, {\Bbb
B}^n),$ where $\Vert Q\Vert_1$ denotes $L^1$-norm of $Q$ in ${\Bbb
R}^n,$ $C_n>0$ is some constant depending only on $n,$ and $r_0=d(K,
\partial {\Bbb B}^n).$ }
\end{corollary}

\section{The case when a majorant is integrable over some set of spheres}

\begin{figure}[h]
\centerline{\includegraphics[scale=0.7]{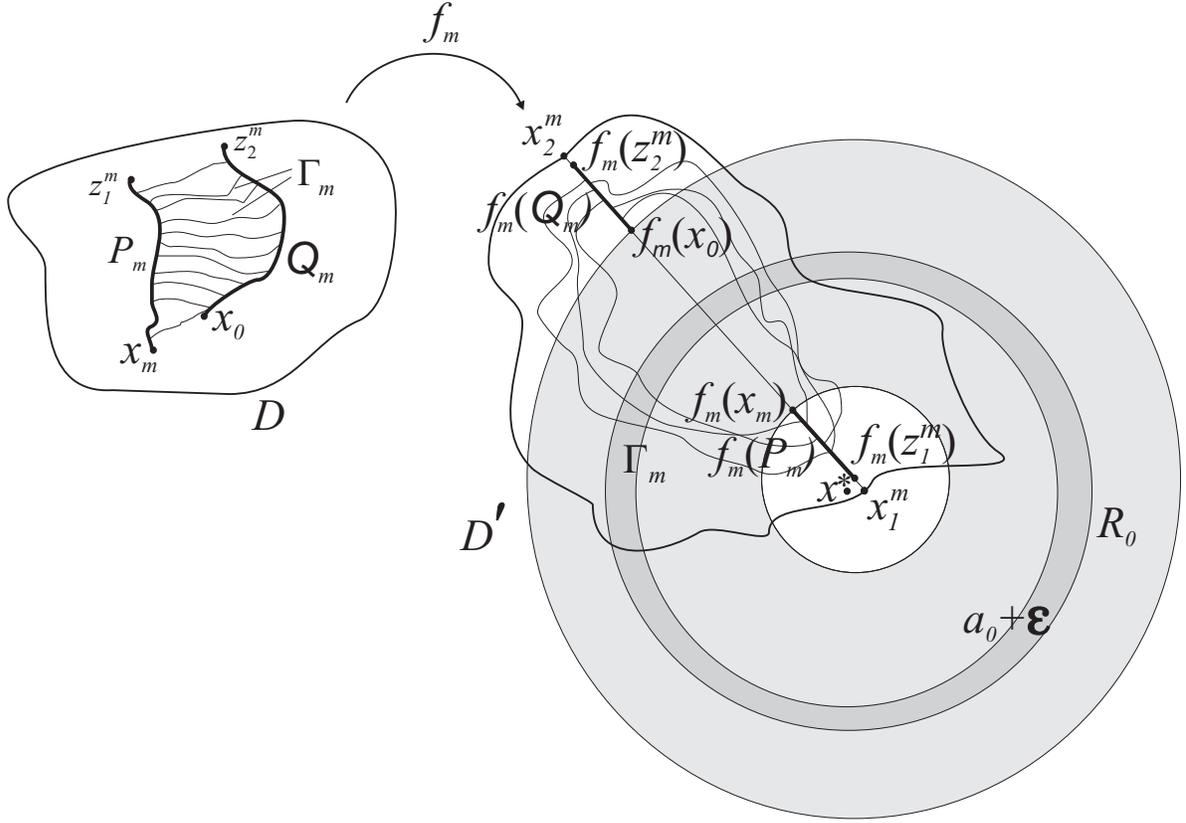}} \caption{To
the proof of Theorem~\ref{th1}}\label{fig2*}
\end{figure}

{\it Proof of Theorem~\ref{th1}.} We prove the theorem~\ref{th1} by
contradiction. Suppose that the conclusion of this theorem does not
hold, that is, the family of maps ${\frak R}_Q(D, D^{\,\prime})$ is
not equicontinuous at some point $x_0\in D.$ Then for each
$m\in{\Bbb N}$ there are $x_m\in D$ and a mapping $f_m\in{\frak
R}_Q(D, D^{\,\prime}),$ such that $|x_m-x_0|<1/m,$ however,
\begin{equation}\label{eq13***}
|f_m(x_m)-f_m(x_0)|\geqslant \varepsilon_0\,.
\end{equation}
Consider a straight line
$$r=r_m(t)=f_m(x_0)+(f_m(x_m)-f_m(x_0))t,\quad-\infty<t<\infty,$$
passing through points $f_m(x_m)$ and $f_m(x_0),$
see~Figure~\ref{fig2*}.
Since the domain $D^{\,\prime}$ is bounded, there are points $x_1^m$
and $x_2^m$ of intersection of the line $r_m$ with the boundary of
the domain $D^{\,\prime},$ see, e.g., \cite[Theorem~1.I.5.46]{Ku}.
Let $\gamma_m^1:[1, c_m)\rightarrow D,$ $1<c_m\leqslant \infty,$ be
a maximal $f_m$-lifting of $r=r(t),$ $t\geqslant 1,$ starting at
$x_m,$ existing by~\cite[Lemma~3.12]{MRV$_2$}. By the same lemma
\begin{equation}\label{eq3}
h(\gamma_m^1(t), \partial D)\rightarrow 0
\end{equation}
as $t\rightarrow c_m-0.$ Similarly, let $\gamma_m^2:(d_m,
0]\rightarrow D,$ $-\infty\leqslant d_m<0$ be a maximal
$f_m$-lifting of $r=r_m(t),$ $t\leqslant 0,$ ending at $x_0,$ which
exists by~\cite[Lemma~3.12]{MRV$_2$}. As in~(\ref{eq3}), we have
that
\begin{equation}\label{eq4}
h(\gamma_m^2(t), \partial D)\rightarrow 0
\end{equation}
as $t\rightarrow d_m+0.$ By~(\ref{eq3}) and (\ref{eq4}), there exist
$1\leqslant t_m^1<c_m$ і $d_m\leqslant t_m^2<0$ such that $h(z_1^m,
\partial D)<1/m$ and $h(z_2^m,
\partial D)<1/m,$ where $\gamma_m^1(t_m^1)=z_1^m$ and $\gamma_m^2(t_m^2)=z_2^m.$
Set $P_m=\gamma_m^1|_{[1, t_m^1]},$ $Q_m=\gamma_m^2|_{[t_m^2, 0]}$
and $\Gamma_m:=\Gamma(P_m, Q_m, D).$ Since the inner points of any
domain are weakly flat (see, e.g., \cite[Lemma~2.2]{SevSkv}), we
obtain that
\begin{equation}\label{eq1B}
M(\Gamma_m)=M(\Gamma(P_m, Q_m, D))\rightarrow\infty\,,\quad
m\rightarrow\infty\,.
\end{equation}
We show that condition~(\ref{eq1B}) leads to a contradiction with
the definition of the class of mappings ${\frak R}_Q(D,
D^{\,\prime})$ in~(\ref{eq2*A}). We denote
$$a_m:=|f_m(z_1^m)-f_m(x_m)|,\, b_m=|f_m(z_1^m)-f_m(x_0)|\,.$$
Obviously, by construction
$$f_m(P_m)\subset B(f_m(z_1^m), a_m)\,,$$
\begin{equation}\label{eq4A}
f_m(Q_m)\subset {\Bbb R}^n\setminus B(f_m(z_1^m), b_m)\,.
\end{equation}
Since the domain $D^{\,\prime}$ is bounded, we may assume that all
three sequences $f_m(z_1^m),$ $f_m(x_0)$ and $f_m(x_m)$ converge to
some elements $\widetilde{x},$ $\widetilde{x_1}$ and
$\widetilde{x_2}$ at $m\rightarrow\infty,$ respectively. Now
$a_m\rightarrow a_0$ and $b_m\rightarrow b_0$ as
$m\rightarrow\infty,$ where $a_0=|\widetilde{x}-\widetilde{x_2}|$
and $b_0=|\widetilde{x}-\widetilde{x_1}|.$

\medskip
Fix
\begin{equation}\label{eq7A}
0<\varepsilon<\varepsilon_0/2\,.
\end{equation}
Let $x^{\,*}\in D$ be such that
$|x^{\,*}-\widetilde{x}|<\varepsilon/3.$ Besides that, let
$m_0\in{\Bbb N}$ be such that
\begin{equation}\label{eq8A}
|f_m(z_1^m)-\widetilde{x}|<\varepsilon/3\,,\quad
a_0-\varepsilon/3<a_m<a_0+\varepsilon/3\quad\forall\,\, m>m_0\,.
\end{equation}

Let $x\in B(f_m(z_1^m), a_m).$ Now, by~(\ref{eq8A}) and the triangle
inequality we obtain that
$$|x-x^{\,*}|\leqslant |x-f_m(z_1^m)|+|f_m(z_1^m)-\widetilde{x}|+
|\widetilde{x}-x^{\,*}|<$$
\begin{equation}\label{eq5}
< a_m+\varepsilon/3+\varepsilon/3< a_0+\varepsilon\,,\quad m>m_0\,.
\end{equation}
It follows from~(\ref{eq5}) that
\begin{equation}\label{eq6A}
B(f_m(z_1^m), a_m)\subset B(x^{\,*}, a_0+\varepsilon)\,, \quad
m>m_0\,.
\end{equation}

Let $R_0$ be such that
\begin{equation}\label{eq6B}
a_0+\varepsilon<R_0<a_0+\varepsilon_0-\varepsilon\,.
\end{equation}
Note that the choice of the number $R_0$ in the formula~(\ref{eq6B})
is possible by the definition of the number $\varepsilon$
in~(\ref{eq7A}). Let $y\in B(x^{\,*}, R_0).$ Again, by the triangle
inequality, and by~(\ref{eq8A}) and~(\ref{eq6B}) we obtain that
$$|y-f_m(z_1^m)|\leqslant
|y-x^{\,*}|+|x^{\,*}-\widetilde{x}|+|\widetilde{x}-f_m(z_1^m)|<$$
\begin{equation}\label{eq5a}
<R_0+2\varepsilon/3<a_0+\varepsilon_0-\varepsilon+2\varepsilon/3=
a_0+\varepsilon_0-\varepsilon/3<a_m+\varepsilon_0\leqslant b_m\,.
\end{equation}

It follows from~(\ref{eq5a}) that
\begin{equation}\label{eq9A}
B(x^{\,*}, R_0)\subset B(f_m(z_1^m),  b_m)\,,\quad m>m_0\,.
\end{equation}
Taking into account relations~(\ref{eq4A}), (\ref{eq6A}) and
(\ref{eq9A}), as well as~\cite[Theorem~1.I.5.46]{Ku}, we obtain that
\begin{equation}\label{eq10}
\Gamma_m>\Gamma_{f_m}(x^{\,*}, a_0+\varepsilon, R_0)\,,\quad
m>m_0\,.
\end{equation}
In turn, from relation~(\ref{eq10}), as well as from the definition
of the class~$\frak{R}_Q(D, D^{\,\prime}),$ it follows that
\begin{equation}\label{eq11}
M(\Gamma_m)\leqslant M(\Gamma_{f_m}(x^{\,*}, a_0+\varepsilon,
R_0))\leqslant \int\limits_{A} Q(y)\cdot \eta^n (|y-x^{\,*}|)\,
dm(y) \,,\quad m>m_0\,,
\end{equation}
where $A=A(x^{\,*}, a_0+\varepsilon, R_0),$ and $\eta$ is an
arbitrary non-negative Lebesgue measurable function satisfying
condition~(\ref{eqA2}) for $r_1=a_0+\varepsilon$ and $r_2=R_0.$
Below we use the standard conventions $a/\infty=0$ for $a\ne\infty$
and $a/0=\infty$ if $a>0$ and $0\cdot\infty=0$ (see e.g.
\cite[3.I]{Sa}). Put now
\begin{equation}\label{eq13}
I=\int\limits_{a_0+\varepsilon}^{R_0}\frac{dt}{tq_{x^{\,*}}^{1/(n-1)}(t)}\,,
\end{equation}
where
\begin{equation}\label{eq12}
q_{y_0}(r)=\frac{1}{\omega_{n-1}r^{n-1}}\int\limits_{S(y_0,
r)}Q(y)\,d\mathcal{H}^{n-1}(y)\,, \end{equation} $\omega_{n-1}$ is
the area of the unit sphere ${\Bbb S}^{n-1}$ in ${\Bbb R}^n,$ and
$q_{x^{\,*}}(t)$ is defined in~(\ref{eq12}) for $y_0:=x^{\,*}.$ By
the hypothesis of the theorem, there is a set $E\subset
[a_0+\varepsilon, R_0]$ of positive linear measure such that the
function $q_{x^{\,*}}(t)$ is finite for $t\in E.$ Thus, $I\ne 0$
in~(\ref{eq13}). In this case, the function
$\eta_0(t)=\frac{1}{Itq_{x^{\,*}}^{1/(n-1)}(t)}$
satisfies~(\ref{eqA2}) for $r_1=a_0+\varepsilon$ and $r_2=R_0.$
Substituting this function on the right side of the
inequality~(\ref{eq11}) and applying Fubini's theorem, we obtain
\begin{equation}\label{eq14}
M(\Gamma_m)\leqslant \frac{\omega_{n-1}}{I^{n-1}}<\infty\,,\quad
m>m_0\,,
\end{equation}
where $\omega_{n-1}$ is the area of the unit sphere ${\Bbb S}^{n-1}$
in ${\Bbb R}^n.$ However, relation~(\ref{eq14})
contradicts~(\ref{eq1B}). The obtained contradiction refutes the
assumption made in~(\ref{eq13***}).~$\Box$

\medskip
To illustrate Theorem~\ref{th1}, we consider the following examples.

\medskip
\begin{example}\label{ex1}
This example is devoted to the study of the case when the result of
Theorem~\ref{th1} (Corollary~\ref{cor1}) can be applied to some
family of mappings, although the corresponding function $Q$ is not
integrable in the domain under consideration. In this case, however,
the so-called Lehto type integral diverges for the function $Q$ (see
e.g.~\cite[(7.50)]{MRSY$_2$}). We consider the following function
$\varphi:(0, 1]\rightarrow {\Bbb R},$ defined as follows:
$$\varphi(t)= \left\{
\begin{array}{rr}
1, & t\in \left(\frac{1}{2k+1}, \frac{1}{2k}\right),\,k=1,2,\ldots\,,\\
\frac{1}{t^n},  &  t\in \left[\frac{1}{2k},\,
\frac{1}{2k-1}\right],\,k=1,2,\ldots\,,
\end{array}
\right. $$
\begin{equation}\label{eq10A}Q(x)=\varphi(|x|)\,,\quad Q:\overline{{\Bbb
B}^n}\setminus\{0\}\rightarrow [0, \infty)\,.
\end{equation}
As usual, put
$$q_0(t)=\frac{1}{\omega_{n-1}t^{n-1}}\int\limits_{S(0, t)}Q(x)\,d\mathcal{H}^{n-1}(x)\,.$$
Using the Fubini theorem, as well as the countable additivity of the
Lebesgue integral, we will have:
$$\int\limits_{{\Bbb B}^n}Q(x)\,dm(x)=\int\limits_0^1\int\limits_{S(0, r)}
Q(x)\,d\mathcal{H}^{n-1}dr=$$
\begin{equation}\label{eq1}
=\omega_{n-1}\int\limits_0^1r^{n-1}\varphi(r)\,dr\geqslant\omega_{n-1}\sum\limits_{k=1}^{\infty}
\int\limits_{1/(2k)}^{1/(2k-1)}\frac{dr}{r}=
\omega_{n-1}\sum\limits_{k=1}^{\infty}\ln\left(\frac{2k}{2k-1}\right)\,.
\end{equation}
Note that the series on the right-hand side of~(\ref{eq1}) diverges.
Indeed, by virtue of the Lagrange's mean value theorem
$\ln\left(\frac{2k}{2k-1}\right)=\ln(2k)-\ln(2k-1)=\frac{1}{\theta(k)}\geqslant
\frac{1}{2k},$ where $\theta(k)\in [2k-1, 2k].$ Since
$\sum\limits_{k=1}^{\infty}\frac{1}{2k}=\infty,$ we obtain that
$\sum\limits_{k=1}^{\infty}\ln\left(\frac{2k}{2k-1}\right)=\infty$
and, consequently, $$\int\limits_{{\Bbb B}^n}Q(x)\,dm(x)=\infty\,.$$
On the other hand,
\begin{equation}\label{eq2A}
\int\limits_0^1\frac{dt}{tq^{1/(n-1)}_{0}(t)}\geqslant\sum\limits_{k=1}^{\infty}
\int\limits_{1/(2k+1)}^{1/(2k)}\frac{dr}{r}=\sum\limits_{k=1}^{\infty}\ln\frac{2k+1}{2k}=\infty\,.
\end{equation}
Define a sequence of mappings $f_m:{\Bbb B}^n\rightarrow {\Bbb R}^n$
by
$$f_m(x)=\frac{x}{|x|}\,\rho_m(|x|)\,,\qquad f_m(0):=0\,,$$
where
$$\rho_m(r)=
\exp\left\{-\int\limits_{r}^1\frac{dt}{tq_{0,
m}^{1/(n-1)}(t)}\right\}\,, \qquad q_{0,
m}(r):=\frac{1}{\omega_{n-1}r^{n-1}}\int\limits_{|x|=r}Q_m(x)\,dS\,,
$$
$$Q_m(x)\quad=\quad \left \{\begin{array}{rr} Q(x) , & \
|x|> 1/m\ ,
\\ 1\ ,  &  |x|\le 1/m\,.
\end{array} \right.$$
Note that each of the mappings $f_m$ is a homeomorphism of the unit
ball ${\Bbb B}^n$ onto itself. Now we show that every $f_m,$
$m=1,2,\ldots,$ satisfies the relation~(\ref{eq8B}) where the
function $Q$ defined by~(\ref{eq10A}). First of all, we note that
each mapping $f_m$ belongs to the class $ACL$; moreover, their norm
and Jacobian are calculated by the relations
\begin{equation*}\label{eq15}\Vert
f_m^{\,\prime}(x)\Vert=\frac{\exp\left\{-\int\limits_{|x|}^1\frac{dt}{tq_{0,
m}^{1/(n-1)}(t)}\right\}}{|x|},\quad |J(x,
f_m)|=\frac{\exp\left\{-n\int \limits_{|x|}^1\frac{dt}{tq_{0,
m}^{1/(n-1)}(t)}\right\}}{|x|^nq_{0, m}^{1/(n-1)}(|x|)}\,,
\end{equation*}
see~\cite[Proof of Theorem~5.2]{IS}. Thus, $f_m\in W_{\rm loc}^{1,
n}({\Bbb B}^n\setminus 0).$ Moreover, the so-called inner dilatation
$K_I(x, f_m)$ of $f_m$ at the point $x$ is calculated as $K_I(x,
f_m)=q_{0, m}(|x|)\leqslant q_0(|x|),$ where $q_{0,
m}(r)=\frac{1}{\omega_{n-1}r^{n-1}}\int\limits_{S(0,
r)}Q_m(x)\,d\mathcal{H}^{n-1}(x)$ (see ibid.). In this case, $f_m$
satisfy relation~(\ref{eq8B}) with $Q=K_I(x, f)=q_0(|x|)$ (see, for
example, \cite[Corollary~8.5 and Theorem~8.6]{MRSY$_2$}.

\medskip
Note that the function $Q,$ extended by zero outside the unit ball,
is integrable over almost all spheres centered at each point $x_0\in
{\Bbb B}^ n,$ since this function is locally bounded in ${\Bbb
B}^n\setminus \{0\}.$ Thus, the maps $g_m=f_m^{\,-1},$
$m=1,2,\ldots$ satisfy all the conditions of Corollary~\ref{cor1}
(more generally, Theorem~\ref{th1}), and by this Corollary the
family of mappings $\{g_m\}_{m=1}^{\infty}$ is equicontinuous in the
unit ball. Moreover, the function $Q$ is not integrable in the unit
ball due to relation~(\ref{eq1}), however, the family
$\{f_m\}_{m=1}^{\infty}$ is equicontinuous in ${\Bbb B}^n$ due to
condition~(\ref{eq2A}) and~\cite[Theorem~7.6]{MRSY$_2$}. However,
the equicontinuity of the family $\{f_m\}_{m=1}^{\infty}$ can be
verified directly; moreover, the equicontinuity of the inverse
family $\{g_m\}_{m=1}^{\infty}$ may be obtained from the fact that
$f_m$ converges to some homeomorphism $f$ as $m\rightarrow\infty$
locally uniformly. In this case, $g_m$ also converges locally
uniformly to some homeomorphism $g=f^{\,-1}$
(see~\cite[Lemma~3.1]{RSS$_1$}).
\end{example}

\medskip
\begin{example}\label{ex2}
Now consider the sequence of functions
$$Q_m(x)\quad=\quad \left \{\begin{array}{rr} \frac{1}{|x|^n} , & \
|x|> 1/m\ ,
\\ 1\ ,  &  |x|\leqslant 1/m\,.
\end{array} \right.$$
In the same way as above, we put
$$f_m(x)=\frac{x}{|x|}\,\rho_m(|x|)\,,\qquad f_m(0):=0\,,$$
where
$$\rho_m(r)=
\exp\left\{-\int\limits_{r}^1\frac{dt}{tq_{0,
m}^{1/(n-1)}(t)}\right\}\,, \qquad q_{0,
m}(r):=\frac{1}{\omega_{n-1}r^{n-1}}\int\limits_{|x|=r}Q_m(x)\,dS\,,
$$
Note that the mappings $f_m,$ as well as inverse mappings
$g_m:=f^{\,-1}_m$ can be written in explicit form, namely,
$$f_m(x)\quad=\quad \left \{\begin{array}{rr}
mx\cdot\exp\{\frac{n-1}{n}\cdot((1/m)^{n/(n-1)}-1\} , & \
|x|\leqslant 1/m\ ,
\\ \frac{x}{|x|}\exp\{\frac{n-1}{n}\cdot(|x|^{n/(n-1)}-1)\}\ ,  &  |x|> 1/m\,,
\end{array} \right.$$
$$g_m(y)= \left \{\begin{array}{rr}
\frac{y}{m}\cdot\exp\{\frac{1-n}{n}((1/m)^{n/(n-1)}-1\} , & \
|y|\leqslant \exp\{\frac{n-1}{n}\cdot((1/m)^{n/(n-1)}-1)\}\ ,
\\ \frac{y}{|y|}\cdot\frac{(n\ln|y|+(n-1))^{(n-1)/n}}{(n-1)^{(n-1)/n}}
\ ,  & |y|> \exp\{\frac{n-1}{n}\cdot((1/m)^{n/(n-1)}-1)\}\,.
\end{array} \right.$$
Reasoning as in Example~\ref{ex1}, it can be shown that the maps
$f_m$ satisfy relation~(\ref{eq8B}) for $D={\Bbb B}^n$ and
$Q(x)=\frac{1}{|x|^n}.$ Note that $\int\limits_{{\Bbb
B}^n}Q(x)\,dm(x)=\infty,$ although in this case we have
$\int\limits_0^1\frac{dt}{tq_0^{1/(n-1)}(t)}<\infty.$ Note that, in
contrast to the mappings from Example~\ref{ex1}, we have a locally
uniform convergence of $g_m$ to $g$ as $m\rightarrow\infty,$ where
$g $ is some continuous mapping. Moreover, the ''direct'' sequence
of mappings of $f_m,$ $m=1,2,\ldots,$ is neither convergent, nor
equicontinuous, in ${\Bbb B}^n.$ This fact about the equicontinuity
of the maps $g_m$ is also the result of Corollary~\ref{cor1}
(Theorem~\ref{th1}), since the function $Q$ has finite mean values
over almost all spheres.
\end{example}

\section{Logarithmic H\"{o}lder continuity of mappings of the unit ball}

\begin{figure}[h]
\centerline{\includegraphics[scale=0.5]{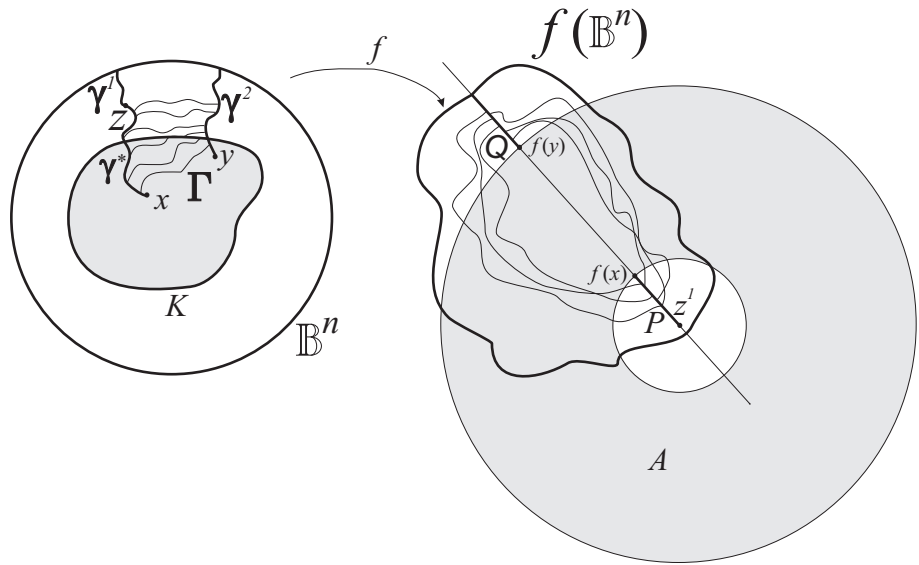}} \caption{To
proof of Theorem~\ref{th1}}\label{fig2}
\end{figure}

{\it Proof of Theorem~\ref{th5}}. We fix $x, y\in K\subset {\Bbb
B}^n$ and $f\in \frak{F}_Q({\Bbb B}^n).$ Put
\begin{equation}\label{eq13A}
|f(x)-f(y)|:=\varepsilon_0\,.
\end{equation}
If $\varepsilon_0=0,$ there is nothing to prove. Let $\varepsilon_0>
0.$
Draw through the points $f(x)$ and $f(y)$ a straight line
$r=r(t)=f(x)+(f(x)-f(y))t,$ $-\infty<t<\infty$ (see
Figure~\ref{fig2}).
Let $\gamma^1:[1, c)\rightarrow {\Bbb B}^n,$ $1<c\leqslant \infty$
be a maximal $f$-lifting of the ray $r=r(t),$ $t\geqslant 1,$
starting at $x,$ which exists by~\cite[Lemma~3.12]{MRV$_2$}. By the
same lemma
\begin{equation}\label{eq3B}
h(\gamma^1(t), \partial {\Bbb B}^n)\rightarrow 0
\end{equation}
as $t\rightarrow c-0.$ Similarly, let $\gamma^2:(d, 0]\rightarrow
{\Bbb B}^n,$ $-\infty\leqslant d<0$ be a maximal $f$-lifting of the
ray $r=r(t),$ $t\leqslant 0,$  with end at a point $y,$ which exists
by~\cite[Lemma~3.12]{MRV$_2$}. Just like in~(\ref{eq3B}) we have
that
\begin{equation*}\label{eq4B}
h(\gamma^2(t), \partial {\Bbb B}^n)\rightarrow 0
\end{equation*}
as $t\rightarrow d+0.$ Let $z=\gamma_1(t_1)$ be some point on the
path $\gamma_1$ located at the distance $r_0/2$ from the boundary of
the unit ball, where $r_0:= d(K, \partial {\Bbb B}^n).$ Put
$\gamma^*(t):=\gamma_1|_{[1, t_1]}.$ Now, by the triangle
inequality, ${\rm diam}\,(|\gamma^*(t)|)\geqslant r_0/2.$
By~\cite[Lemma~4.3]{Vu$_2$}
\begin{equation}\label{eq7D}
M(\Gamma)\geqslant (1/2)\cdot M(\Gamma(|\gamma^*|, |\gamma^2|, {\Bbb
R}^n))\,,
\end{equation}
and on the other hand, by~\cite[Lemma~7.38]{Vu$_1$}
\begin{equation}\label{eq7B}
M(\Gamma(|\gamma^*|, |\gamma^2|, {\Bbb R}^n))\geqslant
c_n\cdot\log\left(1+\frac1m\right)\,,
\end{equation}
where $c_n>0$ is some constant depending only on $n,$
$$m=\frac{{\rm dist}(|\gamma^*|, |\gamma^2|)}{\min\{{\rm diam\,}(|\gamma^*|),
{\rm diam\,}(|\gamma^2|)\}}\,.$$
Note that ${\rm diam}\,(|\gamma^i|)=\sup\limits_{x,
y\in|\gamma^i|}|x-y|\geqslant r_0/2,$ $i=1,2.$ Then
combining~(\ref{eq7D}) and~(\ref{eq7B}), and taking into account
that ${\rm dist}\,(|\gamma^*|, |\gamma^2|)\leqslant |x-x_0|,$ we
obtain that
\begin{equation}\label{eq7C}
M(\Gamma)\geqslant \widetilde{c_n}\cdot \log\left(1+\frac{r_0}{2{\rm
dist}(|\gamma^*|, |\gamma^2|)}\right)\geqslant \widetilde{c_n}\cdot
\log\left(1+\frac{r_0}{2|x-x_0|}\right)\,,
\end{equation}
where $\widetilde{c_n}>0$ is some constant depending only on~$n.$

\medskip
Now let us prove the upper bound for $M(\Gamma).$ Let $z_1:=f(z),$
$P=|f(\gamma^*)|$ and $Q=|f(\gamma^2)|.$ Put
$$A:=A(z^1, \varepsilon^1, \varepsilon^2)=\{x\in {\Bbb R}^n: \varepsilon^1<|x-z^1|<\varepsilon^2\}\,,$$
where $\varepsilon^1:=|f(x)- z^1|,$ $\varepsilon^2:=|f(y)-z^1|.$
We show that
\begin{equation}\label{eq1A}
f(\Gamma)>\Gamma(S(z^1, \varepsilon^1), S(z^1, \varepsilon^2), A)\,.
\end{equation}
Indeed, let $\gamma\in\Gamma.$ Then $f(\gamma)\in f(\Gamma),$
$f(\gamma)=f(\gamma(s)):[0, 1]\rightarrow {\Bbb R}^n,$
$f(\gamma(0))\in P,$ $f(\gamma(1))\in Q$ and $f(\gamma(s))\in
f({\Bbb B}^n)$ for $0<s<1.$ Let $q> 1$ be a number such that
$$z^1=f(y)+(f(x)-f(y))\,q.$$
Since $f(\gamma(0))\in P,$ there is $1\leqslant t\leqslant q$ such
that $f(\gamma(0))=f(y)+(f(x)-f(y))t.$ So,
$$|f(\gamma(0))-z^1|=|(f(x)-f(y))(q-t)|\leqslant $$
\begin{equation}\label{eq2B}\leqslant
|(f(x)-f(y))(q-1)|=|(f(x)-f(y))q+f(y)-f(x))|=\end{equation}
$$=|f(x)-z^1|=\varepsilon^1\,.$$

On the other hand, since $f(\gamma(1))\in Q,$ there is $p\leqslant
0$ such that
$$f(\gamma(1))=f(y)+(f(x)-f(y))p\,.$$ In this case, we obtain that
$$
|f(\gamma(1))-z^1|=|(f(x)-f(y))(q-p)|\geqslant $$
\begin{equation}\label{eq3A}
\geqslant |(f(x)-f(y))q|=|(f(x)-f(y))q +f(y)-f(y)|=
\end{equation}
$$=|f(y)-z^1|=\varepsilon^2\,.$$
Observe that
$$|f(y)- f(x)|+\varepsilon^1=$$
\begin{equation}\label{eq5B}
=|f(y)- f(x)|+|f(x)-z^1|= |z^1-f(y)|=\varepsilon^2\,,
\end{equation}
and, thus, $\varepsilon^1<\varepsilon^2.$ Now, by~(\ref{eq3A}) we
obtain that
\begin{equation}\label{eq4C}
|f(\gamma(1))-z^1|>\varepsilon^1\,.
\end{equation}
It follows from~(\ref{eq2B}) and (\ref{eq4C}) that $|f(\gamma)|\cap
\overline{B(z^1, \varepsilon^1)}\ne\varnothing\ne (f({\Bbb
B}^n)\setminus \overline{B(z^1, \varepsilon^1)})\cap|f(\gamma)|.$ In
this case, by~\cite[Theorem~1.I.5.46]{Ku} there is $t_1\in (0, 1)$
such that~$f(\gamma(t_1))\in S(z^1, \varepsilon^1).$ We may assume
that $f(\gamma(t))\not\in B(z^1, \varepsilon^1)$ for $t\in (t_1,
1).$ Put $\alpha^1:=f(\gamma)|_{[t_1, 1]}.$

\medskip
On the other hand, since $\varepsilon^1<\varepsilon^2$ and
$f(\gamma(t_1))\in S(z^1, \varepsilon^1),$ we obtain that
$|\alpha^1|\cap B(z^1, \varepsilon^2)\ne\varnothing.$
By~(\ref{eq3A}) we obtain that $(f({\Bbb B}^n)\setminus B(z^1,
\varepsilon^2))\cap |\alpha^1|\ne\varnothing.$ Thus,
by~\cite[Theorem~1.I.5.46]{Ku} there is $t_2\in (t_1, 1)$ such that
$\alpha^1(t_2)\in S(z^1, \varepsilon^2).$ We may assume that
$f(\gamma(t))\in B(z^1, \varepsilon^2)$ for $t\in (t_1, t_2).$ Set
$\alpha^2:=\alpha^1|_{[t_1, t_2]}.$ Now, $f(\gamma)>\alpha^2$ and
$\alpha^2\in \Gamma(S(z^1, \varepsilon^1), S(z^1, \varepsilon^2),
A).$ Thus, (\ref{eq1A}) is proved.

\medskip
It follows from~(\ref{eq1A}) that $\Gamma>\Gamma_{f}(z^1,
\varepsilon^1, \varepsilon^2).$ Now, we set
$$\eta(t)= \left\{
\begin{array}{rr}
\frac{1}{\varepsilon_0}, & t\in [\varepsilon^1, \varepsilon^2],\\
0,  &  t\not\in [\varepsilon^1, \varepsilon^2]\,.
\end{array}
\right. $$
Note that $\eta$ satisfies the relation~(\ref{eqA2}) for
$r_1=\varepsilon^1$ and $r_2=\varepsilon^2.$ Indeed, it follows
from~(\ref{eq13A}) and (\ref{eq5B}) that
$$r_1-r_2=\varepsilon^2-\varepsilon^1=|f(y)-z^1|-|f(x)-
z^1|=$$$$=|f(x)-f(y)|=\varepsilon_0\,.$$ Then
$\int\limits_{\varepsilon^1}^{\varepsilon^2}\eta(t)\,dt=(1/\varepsilon_0)\cdot
(\varepsilon^2-\varepsilon^1)\geqslant 1.$ By the
inequality~(\ref{eq1A}), and the relation~(\ref{eq2*A}), applied at
the point~$z^1,$ we obtain that
$$M(\Gamma)\leqslant M(\Gamma_{f}(z^1, \varepsilon^1, \varepsilon^2)) \leqslant$$
\begin{equation}\label{eq14***}
\leqslant \frac{1}{\varepsilon_0^n}\int\limits_{{\Bbb R}^n}
Q(z)\,dm(z)=\frac{\Vert Q\Vert_1}{{|f(x)-f(y)|}^n}\,.
\end{equation}
By~(\ref{eq7C}) and (\ref{eq14***}), we obtain that
$$\widetilde{c_n}\cdot \log\left(1+\frac{r_0}{2|x-y|}\right)\leqslant
\frac{\Vert Q\Vert_1}{{|f(x)-f(y)|}^n}\,.$$
From the latter ratio, the desired inequality~(\ref{eq2C}) follows,
while $C_n:={\widetilde{c_n}}^{-1/n}.$~$\Box$

\section{On logarithmic H\"{o}lder type maps in arbitrary domains}

As we indicated above, the problem of H\"{o}lder continuity for
mappings of arbitrary domains has not yet been resolved.
Nevertheless, under certain not too rough conditions, a result of
this kind can be obtained from the main theorem of the previous
section. In order to formulate and prove it, we carry out the
following notation.

\medskip
For domains $D, D^{\,\prime}\subset {\Bbb R}^n,$ $n\geqslant 2,$ and
a Lebesgue measurable function $Q:{\Bbb R}^n\rightarrow [0,
\infty],$ $Q(x)\equiv 0$ for $x\not\in D,$ we define by
$\frak{A}_Q(D, D^{\,\prime})$ the family of all open discrete
mappings $f:D\rightarrow D^{\,\prime}$ such that
relation~(\ref{eq2*A}) holds at each point $y_0\in D^{\,\prime}.$
The following result holds.

\medskip
\begin{theorem}\label{th6}
{\sl Let $n\geqslant 2,$ and let $Q\in L^1({\Bbb R}^n).$ Suppose
that $K$ is a compact set in $D,$ and beside that, $D^{\,\prime}$ is
bounded. Now, the inequality
\begin{equation*}\label{eq2E}
|f(x)-f(y)|\leqslant\frac{C\cdot (\Vert
Q\Vert_1)^{1/n}}{\log^{1/n}\left(1+\frac{r_0}{2|x-y|}\right)}
\end{equation*}
holds for any $x, y\in K$ and any $f\in \frak{A}_Q(D,
D^{\,\prime}),$ where $\Vert Q\Vert_1$ denotes $L^1$-norm of $Q$ in
${\Bbb R}^n,$ $C=C(n, K)>0$ is some constant depending only on $n$
and $K,$ and $r_0=d(K, \partial D).$ }
\end{theorem}

\medskip
\begin{proof}
First of all, we observe that the equality $\Vert Q\Vert_1= 0$
cannot hold under the conditions of Theorem~\ref{th6}. Indeed, let
$B(x_0, r),$ $r> 0,$ be an arbitrary ball lying strictly inside the
domain $D.$ Let also $\varphi$ be the conformal mapping of the unit
ball ${\Bbb B}^n$ onto the ball $B(x_0, r).$ Applying the
restriction of the mapping $\widetilde{f}:=f|_{B(x_0, r)}$ and
considering the auxiliary map $F:=\widetilde{f}\circ \varphi$,
$F:{\Bbb B}^n\rightarrow D^{\,\prime},$ we conclude that it also
satisfies condition~(\ref{eq2*A}) with the same $Q.$ Thus, by
Theorem~\ref{th5}, $F$ satisfies estimate~(\ref{eq2C}) for any
compact set $K^{\,\prime}\subset {\Bbb B}^n.$ Now, since $F$ is open
and discrete, $\Vert Q\Vert_1\ne 0,$ as required.

\medskip
In view of what has been said, it suffices to establish that the
expression
$$|f(x)-f(y)|\cdot \log^{1/n}\left(1+\frac{r_0}{2|x-y|}\right)$$
is upper bounded for all $x, y\in K.$

\medskip
Suppose the contrary. Then there are $x_m, y_m \in K$ and
$f_m\in\frak{A}_Q(D, D^{\,\prime})$ such that
\begin{equation}\label{eq15A}
|f_m(x_m)-f_m(y_m)|\cdot
\log^{1/n}\left(1+\frac{r_0}{2|x_m-y_m|}\right)\geqslant m\,,\quad
m=1,2,\ldots .
\end{equation}
Since $K$ is compact set, we may assume that both sequences $x_m,
y_m$ are convergent to some points $x_0, y_0 \in K $ as
$m\rightarrow\infty.$

\medskip
Further, two cases are possible: $x_0=y_0$ and $x_0\ne y_0.$ First
let $x_0\ne y_0.$ Then by the triangle inequality
\begin{equation}\label{eq16}
|f_m(x_m)-f_m(y_m)|\leqslant |f_m(x_m)|+|f_m(y_m)|\leqslant C
\end{equation}
for some $C>0$ and any $m\in{\Bbb N},$ because $D^{\,\prime}$ is
bounded by the assumption. It follows from~(\ref{eq16}) that
$$|f_m(x_m)-f_m(y_m)|\cdot
\log^{1/n}\left(1+\frac{r_0}{2|x_m-y_m|}\right)\leqslant$$
\begin{equation}\label{eq13B}
\leqslant C \cdot
\log^{1/n}\left(1+\frac{r_0}{2|x_m-y_m|}\right)\rightarrow C \cdot
\log^{1/n}\left(1+\frac{r_0}{2|x_0-y_0|}\right)<\infty\,, \quad
m\rightarrow\infty\,.
\end{equation}
Relation~(\ref{eq13B}) contradicts assumption~(\ref{eq15A}), so the
case~$x_0\ne y_0$ is impossible.

\medskip
We now consider another case, namely, when~$x_0=y_0.$ Note that in
this case both sequences $x_m$ and $y_m$ belong to $B(x_0, r_0)$ for
any $m>m_0$ and some $m_0\in {\Bbb N},$ while $\overline{B(x_0,
r_0/2)}\subset B(x_0, r_0)\subset D.$ Let $\psi$ be the conformal
mapping of the unit ball ${\Bbb B}^n$ onto the ball $B(x_0, r_0),$
namely, $\psi(x)=xr_0+x_0,$ $x\in {\Bbb B}^n.$ In particular,
$\psi^{\,-1}(B(x_0, r_0/2))=B(0, 1/2).$ Applying the restriction of
the mapping $\widetilde{f}_m:=(f_m)|_{B(x_0, r_0)}$ and considering
the auxiliary maps $F_m:=\widetilde{f}_m\circ \psi,$ $F_m:{\Bbb
B}^n\rightarrow D^{\,\prime},$ we conclude that it also satisfies
condition~(\ref{eq2*A}) with the same $Q.$ Now, by Theorem~\ref{th5}
\begin{equation}\label{eq2F}
|F_m(\psi^{\,-1}(x_m))-F_m(\psi^{\,-1}(y_m))|\leqslant\frac{C\cdot
(\Vert Q\Vert_1)^{1/n}}{\log^{1/n}
\left(1+\frac{1}{4|\psi^{\,-1}(x_m)-\psi^{\,-1}(y_m)|}\right)}\,,
\quad m>m_0\,.
\end{equation}
Since $F_m(\psi^{\,-1}(x_m))=f_m(x_m)$ and
$F_m(\psi^{\,-1}(y_m))=f_m(y_m),$ we may rewrite~(\ref{eq2F}) as
\begin{equation}\label{eq2G}
|f_m(x_m)-f_m(y_m)|\leqslant\frac{C\cdot (\Vert
Q\Vert_1)^{1/n}}{\log^{1/n}\left(1+\frac{1}{4|\psi^{\,-1}(x_m)-\psi^{\,-1}(y_m)|}\right)}\,,
\quad m>m_0\,.
\end{equation}
Note that the maps $\psi^{\,-1}(y)$ are Lipschitz with the Lipschitz
constant $\frac{1}{r_0}.$ In this case, we obtain from
relation~(\ref{eq2G}) that
\begin{equation}\label{eq2H}
|f_m(x_m)-f_m(y_m)|\leqslant\frac{C\cdot (\Vert
Q\Vert_1)^{1/n}}{\log^{1/n}\left(1+\frac{r_0}{4|x_m-y_m|}\right)}\,,
\quad m>m_0\,.
\end{equation}
Finally, we note that
$\log^{1/n}\left(1+\frac{1}{nt}\right)\sim\log^{1/n}\left(1+\frac{1}{kt}\right)$
as $t\rightarrow+0$ for various fixed $k, n> 0,$ which may be
verified using the L'H\^{o}pital rule. It follows that
$$\frac{C\cdot (\Vert
Q\Vert_1)^{1/n}}{\log^{1/n}\left(1+\frac{r_0}{4|x_m-y_m|}\right)}\leqslant
\frac{C_1\cdot (\Vert
Q\Vert_1)^{1/n}}{\log^{1/n}\left(1+\frac{r_0}{2|x_m-y_m|}\right)}$$
for some constant $C_1> 0,$ some $m_1\in {\Bbb N}$ and any
$m>m_1>m_0.$ In this case, it follows from~(\ref{eq2H}) that
\begin{equation}\label{eq2I}
|f_m(x_m)-f_m(y_m)|\leqslant\frac{C_1\cdot (\Vert
Q\Vert_1)^{1/n}}{\log^{1/n}\left(1+\frac{r_0}{2|x_m-y_m|}\right)}\,,
\quad m>m_1\,.
\end{equation}
However, relation~(\ref{eq2I}) contradicts~(\ref{eq15A}). The
obtained contradiction proves the theorem.~$\Box$
\end{proof}

\section{Applications to the Beltrami equation}

Recently, the topic related to the existence of solutions of
degenerate Beltrami differential equations has been actively
developed (see, e.g., \cite{RSY$_1$}, \cite{RSY$_2$}, \cite{GRY},
and \cite{GRSY}). The main results on this topic are compiled in a
relatively recent monograph~\cite{GRSY}, with links to publications
by these and other authors. One of the problems posed in the study
of Beltrami equations is to find the conditions for a complex
coefficient that ensure that their solutions exist. Finding
solutions is usually is carried out in the class of
$ACL$-homeomorphisms, although it is quite correct to consider just
continuous $ACL$-solutions. In this section, we obtain another
result on the existence of solutions of degenerate Beltrami
equation, which is based on the transition to inverse mappings.
Compared to~\cite{RSY$_1$}, \cite{RSY$_2$} and \cite{GRY} , we are
somewhat weakening the conditions on a complex coefficient. The
obtained solution of the equation may not be homeomorphic, but
relative to the previous results, the degree of its smoothness
is~$W_{\rm loc}^{1, p},$ $p>1,$ and therefore somewhat higher.

\medskip We turn now to the definitions.
Let $D$ be a domain in ${\Bbb C}.$ In what follows, a mapping
$f:D\rightarrow{\Bbb C}$ is assumed to be {\it sense-preserving,}
moreover, we assume that $f$ has partial derivatives almost
everywhere. Put $f_{\overline{z}} = \left(f_x + if_y\right)/2$ and
$f_z = \left(f_x - if_y\right)/2.$ The {\it complex dilatation} of
$f$ at $z\in D$ is defined as follows:
$\mu(z)=\mu_f(z)=f_{\overline{z}}/f_z$ for $f_z\ne 0$ and $\mu(z)=0$
otherwise. The {\it maximal dilatation} of $f$ at $z$ is the
following function:
\begin{equation}\label{eq1D}
K_{\mu}(z)=K_{\mu_f}(z)=\quad\frac{1+|\mu (z)|}{1-|\mu\,(z)|}\,.
\end{equation}
Note that the Jacobian of $f$ at $z\in D$ may be calculated
according to the relation
$$J(z,
f)=|f_z|^2-|f_{\overline{z}}|^2\,.$$ Since we assume that the map
$f$ is sense preserving, the Jacobian of this map is positive at all
points of its differentiability. Let ${\Bbb D}=\{z\in {\Bbb C}:
|z|<1\},$ and let $\mu:D\rightarrow {\Bbb D}$ be a Lebesgue
measurable function. Without reference to some mapping $f,$ we
define the {\it maximal dilatation} corresponding to its complex
dilatation~$\mu$ by~(\ref{eq1D}).
It is easy to see that
$$K_{\mu_f}(z)=\frac{|f_z|+|f_{\overline{z}}|}{|f_z|-|f_{\overline{z}}|}$$
whenever partial derivatives of $f$ exist at $z\in D$ and, in
addition, $J(z, f)\ne 0.$
We also define the  {\it inner dilatation of the order $p\geqslant 1
$} of the map $f$ by the relation
\begin{equation}\label{eq17}
K_{I, p}(z,
f)=\frac{{|f_z|}^2-{|f_{\overline{z}}|}^2}{{(|f_z|-|f_{\overline{z}}|)}^p}
\end{equation}
whenever $J(z, f)\ne 0,$ in addition, we set $K_{I, p}(z)=1$
provided that $|f_z|+|f_{\overline{z}}|=0,$ and $K_{I, p}(z)=\infty$
when $J(z, f)=0,$ but $|f_z|+|f_{\overline{z}}|\ne 0.$ Observe that
$K_{I, 2}(z)=K_{\mu}(z).$ Set $\Vert
f^{\,\prime}(z)\Vert=|f_z|+|f_{\overline{z}}|.$ Recall that a
homeomorphism $f:D\rightarrow {\Bbb C}$ is said to be {\it
quasiconformal} if $f\in W_{\rm loc}^{1, 2}(D)$ and, in addition,
$\Vert f^{\,\prime}(z)\Vert^2\leqslant K\cdot |J(z, f)|$ for some
constant $K\geqslant 1.$

\medskip
We will call the {\it Beltrami equation} the differential equation
of the form
\begin{equation}\label{eq2J}
f_{\overline{z}}=\mu(z)\cdot f_z\,,
\end{equation}
where $\mu=\mu(z)$ is a given function with $|\mu(z)|<1$ a.a. Given
$k\geqslant 1$ we set
\begin{equation}\label{eq12:} \mu_k(z)= \left
\{\begin{array}{rr}
 \mu(z),& K_{\mu}(z)\leqslant k,
\\ 0\ , & K_{\mu}(z)> k\,.
\end{array} \right.
\end{equation}
Let $f_k$ be a homeomorphic $ACL$-solution of the equation
$f_{\overline{z}}=\mu_k(z)\cdot f_z,$ which maps the unit disk onto
itself and satisfies the normalization conditions $f_k(0)=0,$
$f_k(1)=1.$ This a solution exists by~\cite[Theorem~3.B.V]{A} or
\cite[Theorem~8.2]{Bo}. Note that the inverse mapping
$g_k=f^{\,-1}_k$ is quasiconformal; in particular, it is
differentiable almost everywhere (see, for example, \cite[Theorems
5.3 and 9.1]{BI}). Let $g_k$ be a inverse mapping to $f_k,$ then its
complex dilation $\mu_{g_k}$ is calculated according to the relation
$\mu_{g_k}(w)=-\mu_k(g_k(w))=-\mu_k(f^{\,-1}_k(w)),$ see e.g.,
\cite[(4).C.I]{A}. In this case, the maximal dilation of $g_k$ is
calculated by the relation
\begin{equation}\label{eq3D}
K_{\mu_{g_k}}(w)=\quad\frac{1+|\mu_k(f^{\,-1}_k(w))|}{1-|\mu_k(f^{\,-1}_k(w))|}\,.
\end{equation}
Accordingly, the inner dilatation of the order $p$ of the map $g_k$
can be calculated according to relation~(\ref{eq17}), namely,
\begin{equation}\label{eq18}
K_{I, p}(w,
g_k)=\frac{{|(g_k)_w|}^2-{|(g_k)_{\overline{w}}|}^2}{{(|(g_k)_w|-|(g_k)_{\overline{w}}|)}^p}\,.
\end{equation}

\medskip
Let $\Omega_n$ be a volume of the unit ball ${\Bbb B}^n$ in ${\Bbb
R}^n.$ Suppose that a function ${\varphi}:D\rightarrow{\Bbb R}$ is
locally integrable in some neighborhood of a point $x_0\in D .$ We
say that ${\varphi}$ has a {\it finite mean oscillation} at $x_0\in
D,$ write $\varphi\in FMO(x_0),$ if the relation
\begin{equation*}\label{eq17:}
{\limsup\limits_{\varepsilon\rightarrow
0}}\,\frac{1}{\Omega_n\varepsilon^n}\int\limits_{B(
x_0,\,\varepsilon)}
|{\varphi}(x)-\overline{{\varphi}}_{\varepsilon}|\ dm(x)\, <\,
\infty
\end{equation*}
holds, where
$\overline{{\varphi}}_{\varepsilon}=\frac{1}{\Omega_n\varepsilon^n}\int\limits_{B(
x_0,\,\varepsilon)} {\varphi}(x)\ dm(x)$ (see, e.g.,
\cite[section~2]{RSY$_2$}).
We say that a function $\varphi$  has a finite mean oscillation in
$D,$ write $\varphi\in FMO(D),$ if $\varphi\in FMO(x_0)$ for any
$x_0\in D.$ The following statement holds.

\medskip
\begin{theorem}\label{th1A}{\sl\,
Let $g_k=f_k^{\,-1}$ and $f_k$ be a homeomorphic $ACL$-solution of
the equation $f_{\overline{z}}=\mu_k(z)\cdot f_z,$ that maps the
unit disk onto itself and satisfies the normalization conditions
$f_k(0)=0,$ $f_k(1)=1$; moreover, let $\mu_k(z)$ be given by the
relation~(\ref{eq12:}). Suppose that the functions $Q:{\Bbb
D}\rightarrow[1, \infty]$ and $\mu:{\Bbb D}\rightarrow {\Bbb D}$ are
Lebesgue measurable, in addition, assume that the following
conditions hold:

\medskip
1) for each $0<r_1<r_2<1$ and $y_0\in {\Bbb D}$  there is a set
$E\subset[r_1, r_2]$ of positive linear Lebesgue measure such that
the function $Q$ is integrable over the circles $S(y_0, r)$ for any
$r\in E;$

\medskip
2) there exist a number $1<p\leqslant 2$ and a constant $M>0$ such
that
\begin{equation}\label{eq10B}
\int\limits_{{\Bbb D}}K_{I, p}(w, g_k)\,dm(w)\leqslant M
\end{equation}
for almost all $w\in {\Bbb D}$ and all $k=1,2,\ldots .$ where $K_{I,
p}(z, g_k)$ is defined in~(\ref{eq18});

\medskip
3) the inequality
\begin{equation}\label{eq10C}
K_{\mu_{g_k}}(w)\leqslant Q(w)
\end{equation}
holds for a.e. $w\in {\Bbb D},$ where $K_{\mu_{g_k}}$ is defined
in~(\ref{eq3D}).

Then the equation~(\ref{eq2J}) has a continuous $W_{\rm loc}^{1,
p}({\Bbb D})$-solution $f$ in ${\Bbb D}.$} \end{theorem}

\medskip
\begin{corollary}\label{cor5}
{\sl\, In particular, the conclusion of Theorem~\ref{th1A} holds if,
in this theorem, the conditions on the function $Q$ are replaced by
the condition $Q\in L^1({\Bbb D}).$ In this case, the solution $f$
of equation~(\ref{eq2J}) can be chosen such that
\begin{equation}\label{eq11A}
|f(x)-f(y)|\leqslant\frac{C\cdot (\Vert
Q\Vert_1)^{1/n}}{\log^{1/n}\left(1+\frac{r_0}{2|x-y|}\right)}
\end{equation}
for any compact set $K\subset {\Bbb D}$ and $x\in x, y\in K,$ where
$\Vert Q\Vert_1$ denotes $L^1$-norm of $Q$ in ${\Bbb D},$ $C>0$ is
some constant, and $r_0=d(K,
\partial {\Bbb D}).$}
\end{corollary}

\medskip
\begin{corollary}\label{cor6}
{\sl If, under the conditions of Corollary~\ref{cor5}, we require in
addition that either $Q(z)\in FMO({\Bbb D}),$ or
\begin{equation}\label{eq5**}
\int\limits_{0}^{\delta(w_0)}\frac{dt}{tq_{w_0}(t)}=\infty
\end{equation}
for any~$w_0\in {\Bbb D}$ and some~$\delta(w_0)>0,$
$q_{w_0}(r)=\frac{1}{2\pi}\int\limits_{0}^{2\pi}Q(w_0+re^{\,i\theta})\,d\theta,$
then $f$ can be chosen as a homeomorphism in ${\Bbb D}.$}
\end{corollary}

\medskip
The proof of Theorem~\ref{th1A}, as well as Corollaries~\ref{cor5}
and~\ref{cor6}, will be given later in the text. Before this, we
formulate and prove the following most important convergence lemma.
On this occasion, see also similar statements related to the study
of equations and mappings with some another conditions, see, for
example, \cite[Ch.~2]{GRSY} and \cite{RSY$_3$}.

\medskip
\begin{lemma}\label{lem1}
{\sl\, Let $\mu:D\rightarrow {\Bbb D}$ be a Lebesgue measurable
function. Suppose that $f_k,$ $k=1,2,\ldots $ is a sequence of
sense-preserving $W_{\rm loc}^{1, 2}(D)$-homeomorphisms of $D$ onto
itself with complex coefficients $\mu_k(z).$ Suppose that $f_k$
converges locally uniformly in $D$ to some mapping $f:D\rightarrow
{\Bbb C}$ as $k\rightarrow\infty,$ and the sequence $\mu_k(z)$
converges to $\mu$ as $k\rightarrow\infty$ for almost all $z\in
{\Bbb D}.$ Suppose also that the inverse mappings $g_k:=f_k^{\,-1}$
belong to the class $W_{\rm loc}^{1, 2}(D)$ and, in addition,
$$\int\limits_{D}K_{I, p}(w, g_k)\,dm(w)\leqslant M$$
for some $M>0,$ each $k=1,2,\ldots ,$ and almost all $w\in D.$ Now
$f\in W_{\rm loc}^{1, p}(D)$ and, in addition, $\mu $ is a complex
characteristic of the map $f,$ in other words,
$f_{\overline{z}}=\mu(z)\cdot f_z$ for almost all $z\in D.$
 }
\end{lemma}

\medskip
\begin{proof}
In general, we will follow the scheme described in the proof
of~\cite[Theorem~3.1]{RSY$_3$}, cf.~\cite[Theorem~2.1]{GRSY} and
\cite[Lemma~III.3.5]{Re}. We denote $\partial f=f_z$ and
$\overline{\partial}f=f_{\overline{z}}.$ Let $C$ be an arbitrary
compact set in $D.$ Since, by assumption, the maps~$g_k=f_k^{\,-1}$
belong to $W_{\rm loc}^{1, 2},$ then $g_k$ possess the Luzin
$N$-property (see, for example, \cite[Corollary~B]{MM}). Now, the
Jacobian $J(z, f)$ is almost everywhere nonzero, see, for example,
\cite[Theorem~1]{Pon}. Since $f_k\in W_{\rm loc}^{1, 2},$ a change
of variables in the integral is true (see
e.g.~\cite[Theorem~3.2.5]{Fe}). In this case, we have that
$$\int\limits_{C}{\Vert f^{\,\prime}_k(z)\Vert}^p\,dm(z)=
\int\limits_C \frac{{\Vert f^{\,\prime}_k(z)\Vert}^p}{J(z,
f_k)}\cdot J(z, f_k)=$$
\begin{equation}\label{eq4D}
=\int\limits_{f_k(C)}K_{I, p}(w, g_k)\, dm(w)\leqslant M<\infty\,.
\end{equation}
It follows from~(\ref{eq4D}) that $f\in W_{\rm loc}^{1, p},$
$\partial f_k$ and $\overline{\partial} f_k$ weakly converge in
$L_{\rm loc}^1(D)$ to $\partial f$ and $\overline{\partial} f,$
respectively (see~\cite[Lemma~2.1]{RSY$_3$} and
\cite[Lemma~III.3.5]{Re}).

\medskip
It remains to show that the map $f$ is a solution of the Beltrami
equation $f_{\overline{z}}=\mu(z)\cdot f_z.$ Put
$\zeta(z)=\overline{\partial} f(z)-\mu(z)\partial f(z)$ and show
that $\zeta(z)=0$ almost everywhere. Let $B$ be an arbitrary disk
lying with its closure in $D.$ By the triangle inequality
\begin{equation}\label{eq9}
\left|\int\limits_B\zeta(z)\,dm(z)\right|\leqslant I_1(k)+I_2(k)\,,
\end{equation}
where
\begin{equation}\label{eq7}
I_1(k)=\left|\int\limits_B(\overline{\partial}f(z)-\overline{\partial}f_k(z))\,dm(z)\right|
\end{equation}
and
\begin{equation}\label{eq8}
I_2(k)=\left|\int\limits_B(\mu(z)\partial f(z)-\mu_k(z)\partial
f_k(z))\,dm(z)\right|\,.
\end{equation}
By proved above, $I_1(k)\rightarrow 0$ as $k\rightarrow\infty.$ It
remains to deal with the expression $I_2(k).$ To do this, note that
$I_2(k)\leqslant I^{\,\prime}_2(k)+I^{\,\prime\prime}_2(k),$ where
$$I^{\,\prime}_2(k)=
\left|\int\limits_B\mu(z)(\partial f(z)-\partial
f_k(z))\,dm(z)\right|$$
and
$$I^{\,\prime\prime}_2(k)=
\left|\int\limits_B(\mu(z)-\mu_k(z))\partial
f_k(z)\,dm(z)\right|\,.$$
Due to the weak convergence of $\partial f_k\rightarrow
\partial f$ in $L^1_{\rm loc}(D)$ as $k\rightarrow\infty,$ we obtain
that $I^{\,\prime}_2(k)\rightarrow 0$ for $k\rightarrow\infty,$
since $\mu\in L^{\infty}(D).$ Moreover, for a given $\varepsilon>0$
there is $\delta=\delta(\varepsilon)>0$ such that
\begin{equation}\label{eq5D}\int\limits_E|\partial f_k(z)|\,dm(z)\leqslant
\int\limits_E|\partial f_k(z)-\partial f(z)|\,dm(z)+
\int\limits_E|\partial f(z)|\,dm(z)<\varepsilon\,,
\end{equation}
whenever $m(E)<\delta,$ $E\subset B,$ and $k$ is sufficiently large.

\medskip
Finally, by Egorov’s theorem (see~\cite[Theorem~III.6.12]{Sa}) for
each $\delta> 0$ there exists a set $S\subset B$ such that
$m(B\setminus S)<\delta$ and $\mu_k(z)\rightarrow \mu(z)$ uniformly
on~$S.$ Then $|\mu_k(z)-\mu(z)|<\varepsilon$ for all $k\geqslant
k_0,$ some $k_0=k_0(\varepsilon)\in {\Bbb N}$ and all $z\in S.$ Now,
by~(\ref{eq5D}) and~(\ref{eq4D}), as well as H\"{o}lder's
inequality,
$$I^{\,\prime\prime}_2(k)\leqslant \varepsilon \int\limits_S|\partial
f_k(z)|\,dm(z)+ 2\int\limits_{B\setminus S}|\partial
f_k(z)|\,dm(z)<$$
\begin{equation}\label{eq6}
<\varepsilon\cdot\left\{\left(\int\limits_D K_{I, p}(w, g_k)\,
dm(w)\right)^{1/p}\cdot (m(B))^{(p-1)/p}+2\right\}\leqslant
\end{equation}
$$\leqslant \varepsilon\cdot\left(M^{1/p}\cdot (m(B))^{(p-1)/p}+2\right)\,.$$
for $k\geqslant k_0.$ From~(\ref{eq9}), (\ref{eq7}), (\ref{eq8}) and
(\ref{eq6}) it follows that $\int\limits_B\zeta(z)\,dm(z)=0$ for all
disks $B,$ compactly embedded in $D.$ Based on the Lebesgue theorem
on differentiation of an indefinite integral
(see~\cite[IV(6.3)]{Sa}), it follows that $\zeta(z)=0$ almost
everywhere in $D.$ The lemma is proved.~$\Box$
\end{proof}

\medskip
{\it Proof of Theorem~\ref{th1A}}. Consider a sequence of
complex-valued functions
\begin{equation}\label{eq12A}
\mu_k(z)= \left \{\begin{array}{rr}
 \mu(z),& K_{\mu}(z)\leqslant k,
\\ 0\ , & K_{\mu}(z)>k,
\end{array} \right.
\end{equation}
where $K_{\mu}(z)$ is defined by~(\ref{eq1D}). Note that
$\mu_k(z)\leqslant\frac{k-1}{k+1}<1,$ therefore, the
equation~(\ref{eq2J}), in which instead of $\mu$ on the right side
we take $\mu:=\mu_k,$ and $\mu_k $ is defined by the
relation~(\ref{eq12A}), has a homeomorphic $W_{\rm loc}^{1, 2}({\Bbb
D})$-solution $f_k:{\Bbb D}\rightarrow {\Bbb D}$ with normalizations
$f_k(0)=0,$ $f_k(1)=1,$ which is $k$-quasiconformal in ${\Bbb D}$
(see~\cite[Theorem~3.B.V]{A} or \cite[Theorem~8.2]{Bo}). By the same
theorem, $f_k$ maps the unit disk onto itself; moreover,
$g_k=f^{\,-1}_k$ are also quasiconformal, in particular, they belong
to the class~$W_{\rm loc}^{1, 2}({\Bbb D}).$
By~\cite[Theorem~6.10]{MRSY$_1$} and by~(\ref{eq10C})
\begin{equation*} \label{eq2*B}
M(g_k(\Gamma))\leqslant \int\limits_{\Bbb
D}K_{\mu_{g_k}}(w)\cdot\rho_*^2 (w) \,dm(w)\leqslant
\int\limits_{\Bbb D}Q(w)\cdot\rho_*^2 (w) \,dm(w)
\end{equation*}
for any $k\in {\Bbb N}$ and any path family $\Gamma$ in ${\Bbb D},$
and each function~$\rho_*\in {\rm adm}\,\Gamma.$ By
Theorem~\ref{th1} the family $\{f_k\}_{k=1}^{\infty}$ is
equicontinuous in~${\Bbb D}.$ Thus, by the Arzela-Ascoli theorem
$f_k$ is a normal family of mappings (see
e.g.~\cite[Theorem~20.4]{Va}), in other words, there is a
subsequence $f_{k_l}$ of $f_k,$ converging locally uniformly in
${\Bbb D}$ to some map $f:{\Bbb D}\rightarrow \overline{{\Bbb D}}.$
Note also that $\mu_k(z)\rightarrow \mu(z)$ as $k\rightarrow \infty$
for almost all $z\in {\Bbb D},$ because $|\mu(z)|<1$ a.e. and,
therefore, $K_{\mu}(z)$ in~(\ref{eq1D}) is finite for almost all
$z\in {\Bbb D}.$ Then by~(\ref{eq10B}) and Lemma~\ref{lem1} the map
$f$ belongs to the class $W_{\rm loc}^{1, p}({\Bbb D})$ and, in
addition, $f$ is a solution of~(\ref{eq2J}).~$\Box$

\medskip
{\it Proof of Corollary~\ref{cor5}} easy follows from
Theorem~\ref{th1A}. Indeed, by the Fubini theorem the condition
$Q\in L^1({\Bbb D})$ implies that the integrals $\int\limits_{S(x_0,
r)\cap D}\,Q(x)\,d\mathcal{H}^1(x)$ are measurable functions by $r$
and finite for a.e. $0<r<\infty$ (see
e.g.~\cite[Theorem~8.1.III]{Sa}). In this case,
conditions~(\ref{eq10B}) and (\ref{eq10C}) are simultaneously
satisfied, where $p=2.$ Thus, the existence of a solution of the
equation~(\ref{eq2J}) and its belonging to the class~$W_{\rm
loc}^{1, p}({\Bbb D})$ follow directly from Theorem~\ref{th1A}.

\medskip
In addition to the above, by Theorem~\ref{th5}
$$|f_k(x)-f_k(y)|\leqslant\frac{C\cdot (\Vert
Q\Vert_1)^{1/2}}{\log^{1/2}\left(1+\frac{r_0}{2|x-y|}\right)}\quad\forall\,\,x,
y\in K$$
where $\Vert Q\Vert_1$ is $L^1$-norm of $Q$ in ${\Bbb D},$ $C$ is
some constant and $r_0:={\rm dist}\,(K,
\partial {\Bbb D}).$ Passing here to the limit as
$k\rightarrow \infty ,$ we obtain the relation~(\ref{eq11A}).
Corollary~\ref{cor5} is proved.~$\Box$

\medskip
{\it Proof of Corollary~\ref{cor6}}. Suppose that~$Q\in FMO({\Bbb
D}),$ or that the relation~(\ref{eq5**}) holds. Then the sequence
$g_k$ forms an equicontinuous family of mappings
(see~\cite[Theorems~6.1 and 6.5]{RS}). Therefore, by the
Arzela-Ascoli theorem $g_k$ is a normal family (see
e.g.~\cite[Theorem~20.4]{Va}), in other words, there is a
subsequence $g_{k_l}$ of $g_k$ converging locally uniformly in
${\Bbb D}$ to some map $g:{\Bbb D}\rightarrow \overline{{\Bbb D}}.$
By the normalization conditions, $g_{k_l}(0)=0$ and $g_{k_l}(1)=1$
for all $l=1,2,\ldots .$ Then, by virtue
of~\cite[Theorem~4.1]{RSS$_1$} the mapping $g$ is a homeomorphism in
${\Bbb D},$ in addition, by~\cite[Lemma ~ 3.1]{RSS$_1$} we also have
that $f_{k_l}\rightarrow f=g^{\,- 1}$ as $l\rightarrow\infty$
locally uniformly in ${\Bbb D}.$ Next, we apply reasoning similar to
those used in the case of an integrable function $Q.$ Since
$\mu_k(z)\rightarrow\mu (z)$ as $k\rightarrow \infty $ and for
almost all $z\in{\Bbb D},$ by Lemma~\ref{lem1} the map $f$ belongs
to the class $W_{\rm loc}^{1, 2}({\Bbb D})$ and, in addition, $f$ is
a solution of~(\ref{eq2J}). Inequality~(\ref{eq11A}) follows from
Corollary~\ref{cor5}.~$\Box$

\medskip
\begin{example}\label{ex3}
Let $p=2,$ let $q\geqslant 1$ be an arbitrary number and let
$0<\alpha<2/q.$ As usual, we use the notation $z=re^{i\theta},$
$r\geqslant 0,$ and $\theta\in [0, 2\pi).$ Put
\begin{equation}\label{eq3A*}\mu(z)= \left
\{\begin{array}{rr}
 e^{2i\theta}\frac{2r-\alpha(2r-1)}{2r+\alpha(2r-1)},& 1/2<|z|<1,
\\ 0\ , & |z|\leqslant 1/2\,.
\end{array} \right.
\end{equation}
Using the ratio
$$\mu_f(z)=\frac{\overline{\partial} f}{\partial f}=e^{2i\theta}\frac{rf_r+if_{\theta}}{rf_r-if_{\theta}}\,,$$
see~(11.129) in \cite{MRSY$_2$}, we obtain that the mapping
\begin{equation}\label{eq4A*}f(z)=\left
\{\begin{array}{rr}
 \frac{z}{|z|}(2|z|-1)^{1/\alpha},& 1/2<|z|<1,
\\ 0\ , & |z|\leqslant 1/2
\end{array} \right.
\end{equation}
is a solution of the equation~$f_{\overline{z}}=\mu(z)\cdot f_z,$
where $\mu$ is defined by~(\ref{eq3A*}). Note that the existence of
a solution of this equation is ensured by Corollary~\ref{cor5} (for
this, we verify that all conditions of this Corollary are
satisfied). Note that for $\mu$ in~(\ref{eq3A*}), the corresponding
maximal dilatation $K_{\mu}$ is the function
\begin{equation}\label{eq5A}K_{\mu}(z)=\left
\{\begin{array}{rr}
 \frac{2|z|}{\alpha(2|z|-1)},& 1/2<|z|<1,
\\ 1\ , & |z|\leqslant 1/2
\end{array} \right.\,.
\end{equation}
Let $k>1/\alpha.$ Observe that $K_{\mu}(z)\leqslant k$ for
$|z|\geqslant \frac{1}{2}\cdot \frac{k\alpha}{k\alpha-1}$ and
$K_{\mu}(z)>k$ otherwise. As above, we set
\begin{equation*}\label{eq5E}\mu_k(z)= \left
\{\begin{array}{rr}
 \mu(z),& K_{\mu}(z)\leqslant k,
\\ 0\ , & K_{\mu}(z)> k\,.
\end{array} \right.
\end{equation*}
Observe that mappings
\begin{equation*}\label{eq6C}f_k(z)=\left
\{\begin{array}{rr}
 \frac{z}{|z|}(2|z|-1)^{1/\alpha},& \frac{1}{2}\cdot \frac{k\alpha}{k\alpha-1}<|z|<1,
\\ \frac{z}{\frac{1}{2}\left(\frac{k\alpha}{k\alpha-1}\right)
}\cdot{\left(\frac{1}{k\alpha-1}\right)}^{1/\alpha}\ , &
|z|\leqslant \frac{k\alpha}{k\alpha-1}
\end{array} \right.\,,
\end{equation*}
are solutions of the equation $f_{\overline{z}}=\mu_k(z)\cdot f_z,$
beside that, the inverse mappings $g_k(y)=f_k^{\,-1}(y)$ are
calculated by the relations
\begin{equation}\label{eq7E}g_k(y)=\left
\{\begin{array}{rr}
 \frac{y(|y|^{\alpha}+1)}{2|y|},& \left(\frac{k\alpha}{k\alpha-1}-1\right)^{1/\alpha}<|y|<1,
\\ \frac{y\cdot\frac{k\alpha}{2(k\alpha-1)}}
{\left(\frac{k\alpha}{k\alpha-1}-1\right)^{1/\alpha}} , &
|y|\leqslant\left(\frac{k\alpha}{k\alpha-1}-1\right)^{1/\alpha}
\end{array} \right.\,.
\end{equation}
It follows from~(\ref{eq5A}) that
\begin{equation}\label{eq7F}K_{\mu_k}(z)=\left
\{\begin{array}{rr}
 \frac{4|z|}{2\alpha(2|z|-1)},& \frac{1}{2}\cdot \frac{k\alpha}{k\alpha-1}<|z|<1,
\\ 1\ , & |z|\leqslant \frac 12\cdot\frac{k\alpha}{k\alpha-1}
\end{array} \right.\,.
\end{equation}
We should check that relation~(\ref{eq10B}) holds for some function
$Q$ that is integrable in ${\Bbb D}.$ For this purpose, we
substitute the maps $g_k$ from~(\ref{eq7E}) into the maximal
dilatation~$K_{\mu_k}$ defined by the equality~(\ref{eq7F}). Then
\begin{equation*}\label{eq8C}K_{\mu_{g_k}}(y)=\left
\{\begin{array}{rr}
 \frac{|y|^{\alpha}+1}{\alpha|y|^{\alpha}},&
 \left(\frac{k\alpha}{k\alpha-1}-1\right)^{1/\alpha}<|y|<1\,,
\\ 1\ , & |y|\leqslant\left(\frac{k\alpha}{k\alpha-1}-1\right)^{1/\alpha}
\end{array} \right.\,.
\end{equation*}
Note that $K_{\mu_{g_k}}(y)\leqslant Q(y):=
\frac{|y|^{\alpha}+1}{\alpha|y|^{\alpha}}$ for all $y\in {\Bbb D}.$
Moreover, the function $Q$ is integrable in ${\Bbb D}$ even in the
degree $q,$ and not only in the degree 1 (see the arguments used in
considering~\cite[Proposition~6.3]{MRSY$_2$}). By the construction
$f_k(0)=0$ and $f_k(1)=1.$ Therefore, all the conditions of
Corollary~\ref{cor5} are satisfied with $p=2$, and the map $f=f(z)$
in~(\ref{eq4A*}) may be considered as the desired solution of the
equation $f_{\overline{z}}=\mu(z)\cdot f_z.$ Moreover, it follows
from the proof of this Corollary that the map $f$ is exactly the
solution of the equation indicated there, since it is a locally
uniform limit of the sequence $f_k.$ Note that the map $f$ is not a
homeomorphic solution, it is also not open and discrete.

\medskip
We show that for a function $\mu$ in~(\ref{eq3A*}) there is no
homeomorphic $W_{\rm loc}^{1, 2}({\Bbb D})$-solution of the Beltrami
equation~(\ref{eq2J}). Indeed, let $g:{\Bbb D}\rightarrow {\Bbb D}$
be such a solution. Due to the Riemann mapping theorem, we may
assume that~$g$ maps the unit disk onto itself. Note that $f$ and
$g$ are locally quasiconformal in $\{1/2<|z|<1\},$ therefore, due to
the uniqueness theorem (see~\cite[Proposition~5.5]{GRSY}),
$g=\varphi\circ f,$ where $\varphi$ is some conformal mapping.
Observe that $\varphi$ is defined in the punctured ball~${\Bbb
D}\setminus\{0\},$ because $f(\{1/2<|z|<1\})={\Bbb
D}\setminus\{0\}.$ Thus, $g\circ f^{\,-1}=\varphi$ and, since
$\varphi$ is conformal in ${\Bbb D}\setminus\{0\},$ $\varphi$ has a
continuous extension to the origin. The last condition cannot be
fulfilled, since $f^{\,-1}(y)= \frac{y(|y|^{\alpha}+1)}{2|y|},$ and
$g$ is some automorphism of the unit disk. This contradiction
disproves the assumption that there exists a homeomorphic $W_{\rm
loc}^{1, 2}({\Bbb D})$-solution $g$ of~(\ref{eq2J}).
\end{example}

\medskip
\begin{example}\label{ex4}
In conclusion, we also give an example of the Beltrami equation, the
existence of an $W_{\rm loc}^{1, p}({\Bbb D})$-solution of which is
ensured by Theorem~\ref{th1A} for some $1<p<2,$ although, at the
same time, Corollaries~\ref{cor5} and~\ref{cor6} are not applicable.
For this purpose, we use the already existing construction of family
of mappings from Example~\ref{ex2}.

\medskip
As usual, we use the notation $z=re^{i\theta},$ $r\geqslant 0,$ and
$\theta\in [0, 2\pi).$ Put
\begin{equation}\label{eq3A**}\mu(z)= \left
\{\begin{array}{rr}
 -e^{2i\theta}\frac{\ln r}{1+\ln r},& e^{-1/2}<|z|<1,
\\ 0\ , & |z|\leqslant e^{-1/2}\,.
\end{array} \right.
\end{equation}
Using the ratio
$$\mu_f(z)=\frac{\overline{\partial} f}
{\partial f}=e^{2i\theta}\frac{rf_r+if_{\theta}}{rf_r-if_{\theta}}\,,$$
see~(11.129) in~\cite{MRSY$_2$}, we obtain that the mapping
\begin{equation}\label{eq4A**}f(z)=\left
\{\begin{array}{rr}
 \frac{z}{|z|}(2\ln|z|+1)^{1/2},& e^{-1/2}<|z|<1,
\\ 0\ , & |z|\leqslant e^{-1/2}
\end{array} \right.
\end{equation}
is a solution of the equation~$f_{\overline{z}}=\mu(z)\cdot f_z,$
where $\mu$ is defined by~(\ref{eq3A**}). Note that the existence of
a solution of this equation is ensured by Theorem~\ref{th1A} (for
this, we verify that all conditions of this Theorem are satisfied).
Note that for $\mu$ in~(\ref{eq3A**}), the corresponding maximal
dilatation $K_{\mu}$ is the function
\begin{equation*}\label{eq5A*}K_{\mu}(z)=\left
\{\begin{array}{rr}
 \frac{1}{1+2\ln r}, & e^{-1/2}<|z|<1,
\\ 1\ , & |z|\leqslant e^{-1/2}
\end{array} \right.\,.
\end{equation*}
Observe that $K_{\mu}(z)\leqslant k$ for $|z|\geqslant
e^{\frac{1-k}{2k}}$ and $K_{\mu}(z)>k$ otherwise. As above, we set
\begin{equation*}\label{eq5F}\mu_k(z)= \left
\{\begin{array}{rr}
 \mu(z),& K_{\mu}(z)\leqslant k,
\\ 0\ , & K_{\mu}(z)> k\,.
\end{array} \right.
\end{equation*}
Observe that mappings
\begin{equation*}\label{eq6C*}
f_k(z)=\left
\{\begin{array}{rr}
 \frac{z}{|z|}(2\ln|z|+1)^{1/2},& e^{\frac{1-k}{2k}}<|z|<1,
\\ ze^{\frac{k-1}{2k}}\cdot k^{-1/2}, &
|z|\leqslant e^{\frac{1-k}{2k}}
\end{array} \right.\,,
\end{equation*}
are solutions of the equation $f_{\overline{z}}=\mu_k(z)\cdot f_z,$
besides that, the inverse mappings $g_k(y)=f_k^{\,-1}(y)$ are
calculated by the relations
\begin{equation}\label{eq7E*}g_k(y)=\left
\{\begin{array}{rr}
 \frac{y}{|y|}\cdot e^{\frac{|y|^2-1}{2}},&
 k^{-1/2}<|y|<1,
\\ ye^{\frac{1-k}{2k}}\cdot k^{1/2} , &
|y|\leqslant k^{-1/2}
\end{array} \right.\,.
\end{equation}
It follows from~(\ref{eq4A**}) that
\begin{equation}\label{eq7F*}K_{\mu_k}(z)=\left
\{\begin{array}{rr}
 \frac{1}{1+2\ln|z|},& e^{\frac{1-k}{2k}}<|z|<1,
\\ 1\ , &|z|\leqslant e^{\frac{1-k}{2k}}
\end{array} \right.\,.
\end{equation}
Now, we substitute the maps $g_k$ from~(\ref{eq7E*}) into the
maximal dilatation~$K_{\mu_k}$ defined by the
equality~(\ref{eq7F*}). We obtain that
\begin{equation*}\label{eq8C*}K_{\mu_{g_k}}(y)=\left
\{\begin{array}{rr}
 \frac{1}{|y|^2},&
 k^{-1/2}<|y|<1,
\\ 1\ , & |y|\leqslant k^{-1/2}
\end{array} \right.\,.
\end{equation*}
We observe that $K_{\mu_{g_k}}(y)$ converges pointwise to
$Q(y)=\frac{1}{|y|^2},$ moreover, by direct calculation we may
verify that the function $Q$ is not integrable in the unit disk
${\Bbb D}.$ It also follows from this that there is no other
function $\varphi$ integrable in the unit disk and such that
$K_{\mu_{g_k}}(y)\leqslant \varphi(y)$ a.e. Indeed, if such a
function existed, then passing here to the limit as
$k\rightarrow\infty $ we obtain that $Q(y)\leqslant \varphi(y),$
which contradicts the condition of non-integrability of
$Q(y)=\frac{1}{|y|^2}$ in ${\Bbb D}.$ On the other hand, the
function $Q$ (extended by zero outside the unit ball) is integrable
over almost all circles $S(x_0, r)$ for any $x_0\in {\Bbb D},$
namely, those that do not pass through the origin.

\medskip
To complete the consideration of the example, we still need to
calculate $K_{I, p}$ in~(\ref{eq18}). For this purpose, we may use
the approach taken when
considering~\cite[Proposition~6.3]{MRSY$_2$}. In the notation of
this Proposition, we obtain that for the mapping $g_k,$
$$\delta_{\tau}=\frac{e^{\frac{|y|^2-1}{2}}}{|y|}\,,
\quad\delta_r=|y|\cdot {e^{\frac{|y|^2-1}{2}}}\,,\quad
k^{-1/2}<|y|<1\,.$$
Thus, $\delta_{\tau}\geqslant \delta_r$ and, consequently,
\begin{equation}\label{eq19}K_{I, p}(y,
g_k)=\frac{e^{\left(\frac{|y|^2-1}{2}\right)(2-p)}}{|y|^p}\,,\quad
k^{-1/2}<|y|<1\,.
\end{equation}
Similarly,
\begin{equation}\label{eq20}
K_{I, p}(y,
g_k)=\left(e^{\frac{1-k}{2k}}k^{\frac{1}{2}}\right)^{2-p}\,,\quad
|y|\leqslant k^{-1/2}\,.
\end{equation}
It follows from~(\ref{eq19}) and (\ref{eq20}) that
$$\int\limits_{\Bbb D}K_{I, p}(y,
g_k)\,dm(y)=\int\limits_{|y|\leqslant k^{-1/2}}K_{I, p}(y,
g_k)\,dm(y)+ \int\limits_{k^{-1/2}<|y|<1}K_{I, p}(y, g_k)\,dm(y)=$$
$$=\pi\cdot k^{\,-1}\cdot \left(e^{\frac{1-k}{2k}}k^{\frac{1}{2}}\right)^{2-p}+
\int\limits_{k^{-1/2}<|y|<1}\frac{e^{\left(\frac{|y|^2-1}{2}\right)(2-p)}}{|y|^p}
\,dm(y)\leqslant$$
$$\leqslant\pi\cdot k^{\,-1+\frac{2-p}{2}}+
\int\limits_{0<|y|<1}\frac{dm(y)}{|y|^p}= \pi k^{\frac{-p}{2}}+ 2\pi
\int\limits_{0}^1 r^{1-p}\,dr\leqslant
\pi+2\pi(2-p)^{\,-1}<\infty\,.$$
Thus, (\ref{eq10B}) is fulfilled with $M=\pi+2\pi(2-p)^{\,-1}.$
\end{example}

\medskip
\medskip
{\bf \noindent Evgeny Sevost'yanov} \\
{\bf 1.} Zhytomyr Ivan Franko State University,  \\
40 Bol'shaya Berdichevskaya Str., 10 008  Zhytomyr, UKRAINE \\
{\bf 2.} Institute of Applied Mathematics and Mechanics\\
of NAS of Ukraine, \\
1 Dobrovol'skogo Str., 84 100 Slavyansk,  UKRAINE\\
esevostyanov2009@gmail.com

\medskip
{\bf \noindent Sergei Skvortsov} \\
Zhytomyr Ivan Franko State University,  \\
40 Bol'shaya Berdichevskaya Str., 10 008  Zhytomyr, UKRAINE \\
serezha.skv@gmail.com

\end{document}